\documentclass[11 pt,final]{article}
\usepackage{amsfonts,latexsym} % symboles
\usepackage{yfonts}
\usepackage{stmaryrd}
\usepackage{color}
\usepackage{overpic}
\usepackage{theorem}
\usepackage{amsmath}
\usepackage{psfrag}
\usepackage{multirow}
\usepackage{rotating}
\usepackage{fancyheadings}
\usepackage{amsfonts}
\usepackage{graphicx}
\usepackage{amssymb}
\usepackage{bm}

\def\mathbi#1{{\boldsymbol{#1}}}
\usepackage{graphicx}
\newcommand{\suff}{eps}
\newcommand{\sufftex}{pstex_t}

\newcommand{\twofig}[4]{\begin{figure}[ht]\begin{center}
    \includegraphics[width=3.5cm]{#1.\suff}\hskip 0.9cm \includegraphics[width=3.5cm]{#2.\suff}
    \end{center}
\caption{#3\label{#4}}\end{figure}}

\newcommand{\threefig}[5]{\begin{figure}[ht]\begin{center}
    \includegraphics[width=3.5cm]{#1.\suff}\hskip 0.4cm \includegraphics[width=3.5cm]{#2.\suff}\hskip 0.4cm \includegraphics[width=3.5cm]{#3.\suff}
    \end{center}
\caption{#4\label{#5}}\end{figure}}

\usepackage[color]{showkeys}
\newcounter{cst}
\newcounter{cexp}
\catcode`\@=11
\def \terml#1{T_{\refstepcounter{cexp}\@bsphack \protected@write\@auxout{}%
           {\string\newlabel{#1}{{\thecexp}{\thepage}}}\thecexp}}
\def \ctel#1{C_{\refstepcounter{cst}\@bsphack \protected@write\@auxout{}%
           {\string\newlabel{#1}{{\thecst}{\thepage}}}\thecst}}
\catcode`\@=12
 \def \cter#1{{C_{\ref{#1}}\, }}
 \def \termr#1{T_{\ref{#1}}}
\newtheorem{definition}{Definition}[section]
\newtheorem{lemma}{Lemma}[section]
\newtheorem{remark}{Remark}[section]
\newtheorem{theorem}{Theorem}[section]

\newenvironment{proof}{\noindent {\sc Proof.} }{$\square$ }

\def\dsp{\displaystyle}

\def\be{\begin{equation}}
\def\ee{\end{equation}}

\def\ba{\begin{array}{lllll}}
\def\ea{\end{array}}

\def\beqsys {\be\ba \left \{ \begin{array}{l}}
\def\eeqsys {\end{array} \right . \ea\ee }

\def\beqsysno {\be\ba \left \{ \begin{array}{l}}
\def\eeqsysno {\end{array} \right . \ea\ee}

\def\bary{{\mathcal B}}

\def\bG{\mathbi{G}}

\def\bpsi{\mathbi{\psi}}

\def\centers{{\cal P}}

\def\cv{K}
\def\cvv{L}

\def\d{\mathrm{d}}

\def\disc{{\mathcal D}}

\def\dcvedge{d_{\cv,\edge}}
\def\dcvedgep{d_{\cv,\edge'}}
\def\dcvvedge{d_{\cvv,\edge}}

\def\div{\mathrm{div}}

\def\dr{\partial}

\def\edge{\sigma}
\def\edgep{{\edge'}}
\def\edges{{\cal E}}
\def\edgesint{{\cal E}_{{\rm int}}}
\def\edgesext{{\cal E}_{{\rm ext}}}
\def\edgeshyb{{\cal H}}

\def\edgescv{{\cal E}_K}
\def\edgescvv{{\cal E}_L}

\def\grad{\nabla}

\def\half{{\frac 1 2}}

\def\lap{\Delta}

\def\matrices{{\cal M}_d(\R)}
\def\mcv{\vert K\vert}
\def\medge{\vert\edge\vert}
\def\medgep{\vert\edge'\vert}

\def\mesh{{\mathcal M}}

\def\ncvedge{\mathbi{n}_{\cv,\edge}}

\def\ncvedgep{\mathbi{n}_{\cv,\edge'}}
\def\nnn{{n \in \N}}
\def\N{\mathbb{N}}

\def\normedeu(#1){\|#1\|_{L^2(\Omega)}}
\def\norm#1#2{\Vert #1 \Vert_{#2}}

\def\O{\Omega}

\def\phi{\varphi}

\def\points{{\cal P}}

\def\R{\mathbb{R}}
\def\refe#1{{\rm(\ref{#1})}}

\def\size{\mathrm{size}}

\def\tends{\to}

\def\x{\mathbi{x}}
\def\xcv{{\mathbi{x}}_\cv}
\def\xcvv{{\mathbi{x}}_\cvv}
\def\xedge{{\mathbi{x}}_\edge}

\def\xedgep{{\mathbi{x}}_{\edge'}}
\def\XD{X_\disc}
\def\XD0{X_{\disc,0}}

\def\y{\mathbi{y}}

\begin{document}

\begin{center}

{\bf \Large Discretisation of   heterogeneous and anisotropic diffusion problems \\ 
on general nonconforming meshes \\ 
SUSHI: a scheme using stabilisation and hybrid interfaces \footnote{This work was supported by Groupement MOMAS, CNRS/PACEN}}

\vspace{1cm}

{R. Eymard\footnote{Universit\'e Paris-Est, France, Robert.Eymard@univ-mlv.fr}, T. Gallou\"et\footnote{Universit\'e de Provence, France, Thierry.Gallouet@cmi.univ-mrs.fr} and  R. Herbin\footnote{Universit\'e de Provence, France,
Raphaele.Herbin@cmi.univ-mrs.fr}}

\vspace{.3cm}
%{\small \today}
\end{center}

{ \small {\bf Abstract:} 
A symmetric discretisation scheme for heterogeneous anisotropic diffusion problems on general meshes is developed and studied. 
The unknowns of this scheme are the values at the centre of the control volumes and at some internal interfaces which may for instance be chosen at the diffusion tensor discontinuities. 
The scheme is therefore completely cell-centred if no edge unknown is kept.
It is shown to be accurate on several numerical examples.  
Convergence of the approximate solution to the continuous solution is proved for general (possibly discontinuous) tensors, general (possibly nonconforming) meshes, and with no regularity assumption on the solution. 
An error estimate is then deduced under suitable regularity assumptions on the solution. 

\smallskip

{\bf Keywords : } {Heterogeneous anisotropic diffusion, nonconforming grids, finite volume schemes}
}
\section{Introduction}
Anisotropic heterogeneous diffusion problems arise in a wide range of scientific fields such as hydrogeology, oil reservoir simulation, plasma physics, semiconductor modelling, biology, etc.. 
When implementing  numerical methods for this kind of problem, one needs to  find an approximation of $u$, weak solution to the following equation:
\be 
 -\div (\Lambda(\x)\grad u) = f \hbox{ in } \Omega, \label{eq1}
\ee with boundary condition 
\be 
u=0 \hbox{ on } \dr \Omega \label{cl0}, 
\ee
where we denote by  $\dr \O =  \overline{\O}\setminus\O$ the boundary of the domain $\O$, under the following assumptions:
\be
\O \mbox{ is an open bounded connected polyhedral subset of }\R^d,
\ d\in\N\setminus\{0\},
\label{hypomega}\ee
\be
\Lambda \hbox{ is a measurable function from   } \Omega  \hbox{ to } \matrices,
\label{hyplambda}\ee  
where we denote by $\matrices $ the set of $d\times d $ matrices,  such that for a.e. $\x \in \Omega$, $\Lambda(\x)$ is symmetric, and such that the set of its eigenvalues is included in $ [\underline{\lambda},\overline{\lambda}]$, with $\underline{\lambda}$ and $\overline{\lambda} \in \R$ satisfying $0 < \underline{\lambda}\le \overline{\lambda}$, and
\be
f \in L^2(\O).
\label{hypfg}
\ee
Under these hypotheses, the weak solution of \refe{eq1}--\refe{cl0} is the unique function $u$ satisfying:
\be\left\{\ba
u \in H^1_0(\O),\\
\dsp \int_\O \Lambda(\x)\grad u(\x)\cdot\grad v(\x) \d\x = \int_\O f(\x)  v(\x) \d\x \qquad \forall v \in H^1_0(\O).
\ea\right.\label{ellgenf}\ee

Usual discretisation schemes for Problem \refe{ellgenf} include finite difference, finite element or finite volume methods.
Finite volume methods are actually very  popular in oilreservoir  engineering, a probable reason being that complex coupled physical phenomena may be discretised on the same grids.
The well-known five-point scheme on rectangles (see e.g. \cite{pat-80-num}) and four-point scheme on triangles \cite{her-95-err} are not easily adapted to heterogeneous anisotropic diffusion operators  \cite{her-96-fin}.
A scheme with an enlarged  stencil, which handles anisotropy on meshes satisfying an orthogonality property, was proposed and analysed in \cite{eym-06-cel}. 
Another problem that has to be faced in several fields of applications (such as hydrogeology and oil reservoir engineering) is the fact that the discretisation meshes are imposed by engineering and computing considerations; therefore, we have to deal with distorted and possibly nonconforming meshes. 

A huge literature exists in the engineering setting, so we shall not try to be exhaustive. 
Let us nevertheless mention the finite volume schemes using the   well-known  multipoint flux approximation \cite{aav-96-dis,aav-98-dis,aav-98-ani}. 
These schemes involve the reconstruction of the gradient in order to evaluate the fluxes, which is also the case in \cite{cou-99-con,lep-05-mon}.
Among other approaches let us cite \cite{gui-05-num}, which uses a parametrisation technique.
However, even though these schemes perform well in a number of cases, their convergence analysis often seems to remain out of reach, except under additional geometrical conditions \cite{cou-99-con}. 

More recently, finite volume schemes using interface values have been studied. 
In \cite{eym-07-new} we presented a ``hybrid finite volume" (HVF) scheme for any space dimension, which involves edge unknowns in addition to the usual cell unknowns, and in \cite{dro-06-mix}, a ``mixed finite volume" scheme (MFV) was proposed, which involves the fluxes and the values as unknowns.
This is also the case for the mimetic finite difference (MFD) schemes \cite{bre-05-con,bre-05-fam}, which were introduced previously;  in spite of their name, mimetic schemes are very much in the finite volume spirit, since they rely  on both a flux balance equation and on the local conservativity of the numerical fluxes, that are probably the two ``pillars" of the finite volume philosophy;
but then, finite volume schemes are also often called finite difference schemes in the engineering literature because of  the finite difference approximation of the fluxes. 
In fact, a recent benchmark  \cite{her-08-ben} provided sufficient information to suspect that the methods HFV, MFV and MFD indeed coincide at the algebraic level and establishing this is the aim of ongoing work  \cite{dro-09-com}.
Let us mention that the Raviart-Thomas mixed finite element method, which also involves edge  unknowns, was  generalised to handle distorted hexahedral meshes \cite{kuz-05-con}.
These schemes require the fluxes or edge unknowns as additional values (or as sole values after hybridisation), and they may be more expensive than cell-centred schemes, especially in the 3D case. 
 
In the two-dimensional case, we also mention \cite{ber-07-ver}, which discusses a scheme based on vertex reconstructions, and the family of double mesh schemes \cite{her-03-app, dom-05-fin, boy-08-fin}. 
The generalisation of this type of scheme to 3D is the subject of ongoing work. 

The scheme that we present here is designed on very general polygonal, possibly non-convex and nonconforming meshes, with the following two priorities in mind: 
\begin{itemize}
\item For cost reasons and data structure issues, we wish to obtain a symmetric scheme which is as close as possible to a cell-centred scheme, that is to a scheme involving one unknown per control volume (or grid cell). 
\item  For accuracy reasons, we require the local conservativity of the numerical flux to hold at the interfaces between highly heterogeneous media.
\end{itemize}
In \cite{eym-07-col}, we introduced a cell-centred scheme for the approximation of the Laplace operator on nonconforming grids in the framework of the incompressible Navier-Stokes equations \cite{eym-07-col} and which may be viewed as a low order nonconforming Galerkin approximation.
The scheme (called ``SUCCES"  in \cite{age-08-sym}) was also implemented for anisotropic and heterogeneous problems on general meshes, and was shown to be highly competitive for oil reservoir simulation  in comparison with other  well-known  schemes such as the multiple point flux approximation schemes.
%However, we show in Section \ref{sec.numres} that in the case of highly heterogeneous cases such as 
It is cheaper than the above mentioned  hybrid type schemes (HVF, MFV and MFD) because it is based on cell unknowns only.
However, it is not as accurate as the hybrid schemes for strongly heterogeneous problems, very likely because of the weaker approximation of the normal fluxes at the heterogeneous interfaces.
%In \cite{eym-07-new} we presented a ``hybrid finite volume" scheme for any space dimension, which involves edges unknowns in addition to the usual cell unknowns. 
%This is also the case for the mimetic finite difference schemes \cite{bre-05-con}.
%Along the same line of thought, a ``mixed finite volume" scheme was proposed in \cite{dro-06-mix}.
%These schemes perform quite well but seem rather expensive at first glance, because of the edge unknowns and equations.
In the present work, we construct a discretisation scheme (SUSHI) for any kind of polyhedral mesh, which incorporates the best properties of the cell-centred (SUCCES) and hybrid (HFV) schemes: 
unknowns on the edges are only introduced when needed, for instance when there is strong medium heterogeneity at these edges. 
If the set of edge unknowns is empty, then SUSHI reduces to the above mentioned cell-centred scheme; if unknowns are associated to all internal edges, then SUSHI is the hybrid scheme HFV.

\medskip

The outline of this paper is as follows. In Section \ref{fondam}, we present the guidelines which led us in the construction of convergent schemes on general nonconforming meshes. 
The practical properties of the resulting schemes are shown through numerical examples in Section \ref{secresnum}. 
Then the mathematical analysis of convergence and error estimation are performed in Section \ref{cvstudy}. 
This analysis is based on some discrete functional analytic tools, such as discrete Sobolev inequalities, which are provided in Section \ref{secsobdis}. 
Conclusions and perspectives are discussed in Section~\ref{conc}.

\section{Fundamentals for a class of nonconforming schemes}\label{fondam}
 
Let us first present the desired properties which have led us to the design of the schemes under study:

\begin{enumerate}
\item[(P1)] The schemes must apply on any type of grid: conforming or nonconforming, 2D and 3D (or more, see for instance the frameworks of kinetic formulations or financial mathematics), consisting of control volumes which are only assumed to be polyhedral (the boundary of each control volume is a finite union of subsets of hyperplanes).
\item[(P2)] The matrices of the linear systems  generated  are expected to be sparse, symmetric and positive definite. 
\item[(P3)] We wish to be able to prove the convergence of the family of  discrete solutions to  the solution of the continuous problem as the mesh size tends to 0, and of the family of associate gradients to the  gradient of the solution, with no regularity assumption on the solution of the continuous problem, and to derive error estimates when the analytic solution is regular enough. 
\end{enumerate}

 In order to describe the schemes we now introduce some notations for the space discretisation.

\begin{definition}[Space discretisation]\label{adisc}
Let $\O$ be a polyhedral open bounded connected subset of $\R^d$, with $d\in\N\setminus\{0\}$, and $\dr \O =  \overline{\O}\setminus\O$ its boundary.
A discretisation of $\O$, denoted by $\disc$, is defined as the triplet $\disc=(\mesh,\edges,\centers)$, where:
\begin{enumerate}

\item $\mesh$ is a finite family of nonempty connected open disjoint subsets of $\O$ (the ``control volumes'') such that $\overline{\O}= \dsp{\cup_{K \in \mesh} \overline{K}}$.
For any $K\in\mesh$, let $\dr K  = \overline{K}\setminus K$ be the boundary of $K$; let $\mcv>0$ denote the measure of $K$ and let $h_K$ denote the diameter of $K$.

\item $\edges$ is a finite family of disjoint subsets of $\overline{\O}$ (the ``edges'' of the mesh), such that, for all $\edge\in\edges$,  $\edge$ is a nonempty open subset of a hyperplane of $\R^d$, whose $(d-1)$-dimensional measure $\medge$ is strictly positive.
We also assume that, for  all $K \in \mesh$, there exists  a subset $\edgescv$ of $\edges$
such that $\dr K  = \dsp{\cup_{\edge \in \edgescv}}\overline{\edge} $. 
For any $\edge \in \edges$, we denote by $\mesh_\edge = \{K\in\mesh, \edge\in\edgescv\}$.
We then assume that, for all $\edge\in\edges$, 
either  $\mesh_\edge$ has exactly one element and then $\edge\subset \dr\O$ (the set of these interfaces, called
boundary interfaces, is denoted by $\edgesext$) or  $\mesh_\edge$ has exactly two elements
(the set of these interfaces, called interior interfaces,  is denoted by $\edgesint$).
For all  $\edge\in\edges$, we denote by $\xedge$ the barycentre of $\edge$. 
For all $K \in \mesh$ and $\edge \in \edgescv$, we denote  by $\ncvedge$ the unit vector normal to $\edge$ outward to $K$.

\item $\centers$ is a family of points of $\O$ indexed by $\mesh$, denoted by $\centers = (\xcv)_{K \in \mesh}$, 
such that for all  $K\in\mesh$,  $\xcv\in K$ and $K$ is assumed to be $\xcv$-star-shaped,  which means that for all $\x\in K$, the inclusion $[\x_K,\x]\subset K$ holds. 
Denoting by $d_{K,\edge}$ the Euclidean distance between $\xcv$ and the hyperplane including $\edge$, one assumes
that  $d_{K,\edge}>0$.  
We then denote by $D_{K,\edge}$ the cone with vertex $\xcv$ and basis $\edge$.

\end{enumerate}
 
\end{definition}

\begin{remark}
The above definition applies to a large variety of meshes. Note that no hypothesis is made on the convexity of the control volumes; 
in fact, generalised hexahedra, {\it i.e.} with faces which may be composed of several planar sub-faces may be used.   
Often encountered in subsurface flow simulations, such   hexahedra  may have up to 12 faces (resp. 24 faces) if each non planar face is composed of two triangles (resp. four triangles), but only 6 neighbouring control volumes. 
\end{remark}

\subsection{From a ``hybrid" finite volume scheme\ldots}
The idea of the ``hybrid'' schemes (among them one may include the mixed finite elements, the mixed finite volume or the mimetic finite difference schemes) is to find an approximation to the solution of \refe{eq1}--\refe{cl0} by setting up a system of discrete equations for a family of values $((u_K)_{K \in \mesh}, (u_\edge)_{\edge \in \edges})$ in the control volumes and on the interfaces. 
The number of unknowns is therefore card($\mesh$) + card($\edges$). 
Following the idea of the finite volume framework, Equation \refe{eq1} is integrated over each control volume $K\in\mesh$, which formally gives (assuming sufficient regularity on $u$ and $\Lambda$) the following balance equation on the control volume $K$:
\[
\sum_{\edge\in\edgescv} \left(-\int_\edge \Lambda(\x) \nabla u(\x) \cdot \ncvedge \d \gamma(\x) \right) = \int_K f(\x) \d\x.
\]
The flux $-\int_\edge \Lambda(\x) \nabla u(\x) \cdot \ncvedge \d \gamma(\x)$ is approximated by a function $F_{K,\sigma}(u)$ of the values  $((u_K)_{K \in \mesh}$, $(u_\edge)_{\edge \in \edges})$ at the ``centres'' and  at the interfaces of the control volumes (in all practical cases,  $F_{K,\sigma}(u)$ only depends on $u_K$ and all $(u_\edgep)_{\edgep\in\edgescv}$). 
A discrete equation corresponding to \refe{eq1} is then:
\be
\sum_{\edge\in\edgescv} F_{K,\sigma}(u) = \int_K f(\x)\d\x \qquad \forall K\in\mesh.
\label{schvolfincla}
\ee
The values $u_\sigma$ on the interfaces are then introduced so as to allow 
for a consistent approximation of the normal fluxes in the case of an anisotropic operator and a 
general, possibly nonconforming mesh. 
We thus have card$(\edges)$ supplementary unknowns, and need  card$(\edges)$ equations to ensure that the problem is well posed. 
For the boundary faces or edges, these equations are obtained by writing the discrete counterpart of the 
boundary condition \refe{cl0}:
\be 
u_\edge = 0 \qquad \forall \edge\in\edgesext. \label{clvf}
\ee
Following the finite volume ideas, we may write the continuity of the discrete flux for all interior edges, that is to say:
\be
F_{K,\sigma}(u)+F_{L,\sigma}(u) = 0, \hbox{ for } \sigma \in \edgesint \mbox{ such that } \mesh_\edge = \{K,L\}.
\label{consfludis}
\ee
We now have card($\mesh$) + card($\edgesint$) unknowns and equations. 
\begin{remark}
%re j'ai mieux detaille le fait qu'il s'agit de la moyenne harmonique
In the case $\Lambda(\x) = \lambda(\x){\rm Id}$, on meshes satisfying an orthogonality condition as mentioned in the introduction of this paper (this condition states the orthogonality between the line joining the centres of two neighbouring control volumes with their
common interface, see \cite[Definition 9.1 p. 762]{book}), a consistent numerical flux is obtained using the two-point formula 
$F_{K,\sigma}(u) = \lambda_K \medge(u_K - u_\edge)/\dcvedge$, where $\lambda_K$ is the average value for $\lambda$ in $K$. 
%re
Then, writing \refe{consfludis} for all $\sigma \in \edgesint$ such that $\mesh_\edge = \{K,L\}$, we obtain  $u_\edge$ as a linear combination of  $u_K$ and $u_L$. 
Plugging this expression into \refe{schvolfincla}, we get a scheme with card$(\mesh)$  equations and  card$(\mesh)$ unknowns (see \cite[Section 11.1 pp. 815-820]{book} for more details). 
In the case of a rectangular (resp. triangular) mesh, this is the well-known five points (resp. four points) scheme with harmonic averages of the diffusion.
\end{remark}

With a proper choice of the expression $F_{K,\sigma}(u)$,  which we shall introduce below, this scheme, first introduced in \cite{eym-07-new}, is quite efficient  for the simulation of fluid flow in heterogeneous media (where harmonic averages for $\Lambda$ are preferred to arithmetic averages \cite{azi-79-pet}) and may be shown to converge. 
This scheme does have one drawback: 
since the number of unknowns is the sum of the number of control volumes and of interior interfaces, the resulting scheme is quite expensive (although it is sometimes possible to algebraically eliminate the values at the control volumes, as in the mixed hybrid finite element method, see \cite[pp. 178-181]{bre-91-mix}).  
%Note that preconditioning with 2-point fluxes can be considered as an efficient method for solving the linear systems \cite{refprecond}.
\begin{remark}
Note that in the case of regular conforming simplices (triangles in 2D, tetrahedra in 3D), there is an algebraic possibility to express the unknowns $(u_\edge)_{\edge \in \edges}$ as local affine combinations of the values $(u_K)_{K \in \mesh}$ and therefore to eliminate them \cite{voh-00-equ}. 
The idea is to remark that the linear system constituted by the equations  \refe{schvolfincla} for  all  $K\in \mesh_S$, where $\mesh_S$ is the set  of all simplices sharing the same interior vertex $S$, and \refe{consfludis} for all the interior edges  such that $\mesh_\edge\subset\mesh_S$, presents as many equations as unknowns $u_\edge$, for  $\edge\in \cup_{K\in \mesh_S} \edgescv$. Indeed, the number of edges in $\cup_{K\in \mesh_S} \edgescv$ such that $\mesh_\edge\not\subset\mesh_S$ is equal to the number of control volumes in  $\mesh_S$. 
Unfortunately, there is at this time no general result on the invertibility or the symmetry of the matrix of this system, and this method does not apply to other types of meshes than simplicial meshes.
\end{remark}

In order to reduce the computational cost of the scheme, we developed in \cite{eym-07-col} an idea which is in fact close to the finite element philosophy since we express the finite volume scheme in a weak form; to this end, let us first define the sets $X_{\disc}$ and $X_{\disc,0}$ where the discrete unknowns lie, that is to say: 
\begin{eqnarray}
 X_{\disc}=\{v = ((v_K)_{K \in \mesh}, (v_\edge)_{\edge \in \edges}), v_K \in \R, v_\edge \in \R\},\label{defX} \\
 X_{\disc,0}=\{v \in X_\disc \mbox{ such that } v_\edge = 0 \qquad \forall \edge \in \edgesext\}.
\label{defX0}
\end{eqnarray}
Multiplying, for any  $v\in X_{\disc,0}$, Equation \refe{schvolfincla} by the value $v_K$ of $v$ on the control volume $K$ and summing over $K\in\mesh$ leads to:
\[
\sum_{K\in\mesh} v_K \sum_{\edge\in\edgescv} F_{K,\sigma}(u) = \sum_{K\in\mesh} v_K \int_K f(\x)\d\x.
\]
Using \refe{consfludis}, we get the following discrete weak formulation:
\begin{equation}
\label{fvweak}
\left\{\begin{array}{l}
\mbox{Find } u \in \XD0 \mbox{ such that:}\\
\dsp \langle u,v \rangle_F   =\sum_{K\in\mesh} v_K \int_K f(\x)\d\x, \ \mbox{ for all } v \in \XD0,
\end{array}\right.
\end{equation}
with 
\begin{equation}
\langle u,v \rangle_F = \sum_{K\in\mesh} \sum_{\edge\in\edgescv} F_{K,\sigma}(u) (v_K - v_\edge). 
\label{defbilF}
\end{equation} 
Note that choosing $v \in \XD0$ such that  $v_K = 1$, $v_L = 0$ for any $L \in \mesh, L\not = K$ and $v_\edge = 0$ for any $\edge \in \edges$ yields \refe{schvolfincla}. Similarly, choosing $v \in \XD0$ such that  $v_K = 0$ for any $K \in \mesh$, and $v_\edge = 1$ and $v_\tau = 0$ for any $\tau \in \edges, \tau \not = \sigma$ leads to \refe{consfludis}. 
Therefore the hybrid finite volume scheme \refe{schvolfincla}--\refe{consfludis} is equivalent to the discrete weak formulation \refe{fvweak}.

\subsection{\ldots to a nonconforming finite element scheme\ldots} \label{sec-succes}

We may then choose to use the weak discrete form \refe{defbilF} as an approximation of the bilinear form $a(\cdot,\cdot)$, but with a  space of dimension smaller than that of $\XD0$. 
% re j'ai repris ici, car on y introduit un hmesh non conforme avec la suite
% 
%Let us write for instance \refe{fvweak} replacing $\XD0$ by $H_\mesh = \{(u_K)_{K \in \mesh}, u_K \in \R\}$. 
%The new formulation therefore amounts to card$(\mesh)$ discrete equations, 
%which is fine so long as we have
%
This can be achieved by expressing 
the value of $u$ on any interior interface $\edge \in \edgesint$ as a consistent barycentric combination of the values $u_K$:
\be
u_\edge = \sum_{K\in\mesh}\beta_\edge^K u_K,
\label{ecrbar}
\ee
where $(\beta_\edge^K)_{\stackrel{K\in\mesh}{\edge\in\edgesint}}$ is a family of real numbers,  with  $\beta_\edge^K \not=  0$ only for some control volumes $K$ close to $\edge$, and such that
\be
 \sum_{K\in\mesh}\beta_\edge^K = 1\hbox{ and }\xedge = \sum_{K\in\mesh}\beta_\edge^K \xcv \qquad\forall \edge \in \edgesint.
\label{propbary}\ee
This ensures that if $\phi$  is a regular function, then $\phi_\sigma = \sum_{K\in\mesh}\beta_\edge^K \phi(\xcv) $ is a consistent approximation of $\phi(\xedge)$ for $\edge \in \edgesint$. 
We recall that the values $u_\edge, \edge \in \edgesext$ are set to 0 in order to respect the boundary conditions \eqref{cl0}.
Hence the new scheme reads: 
\begin{equation}
\left\{
\begin{array}{l}
%%%%%%%%%  \mbox{Find }u \in H_\mesh \mbox{( that is  } (u_K)_{K \in \mesh}),  \mbox{ such that:} \\
\mbox{Find } u \in \XD0  \mbox{ such that } 
\dsp u_\edge = \sum_{K\in\mesh} \beta_\edge^K u_K\qquad \forall \sigma \in \edgesint, \mbox{ and }\\
\dsp \langle u,v \rangle_F   = \sum_{K\in\mesh} v_K \int_K f(\x)\d\x,\ \mbox{ for all } v \in \XD0 \mbox{ with } v_\edge = \sum_{K\in\mesh} \beta_\edge^K v_K\qquad \forall \sigma \in  \edgesint.
\end{array}\right.
\label{succes}
\end{equation}

This method has been shown in \cite{eym-07-col} to be efficient in the case of a problem where $\Lambda = {\rm Id}$ (for the approximation of the viscous terms in the Navier-Stokes problem).
With an appropriate choice for the expression of the numerical flux, it also yields conservativity in a certain sense (more on this below), but no longer to the classical (in the finite volume framework) equation \refe{consfludis}: indeed, since the degrees of freedom on the edges are no longer present, one may not use $v_\sigma = 1$ to recover \refe{consfludis}. Note also that taking $v_K= 1$ does not yield \refe{schvolfincla}. This scheme has been implemented for the discretisation of the diffusive term in the incompressible Navier Stokes equations on general two- or three-dimensional grids, and gives excellent results \cite{che-08-colf,che-09-col}.
Unfortunately, because of poor approximation of the local flux at strongly heterogeneous interfaces, this approach is not sufficient to provide accurate results  for some types of flows in heterogeneous media, as we shall show in Section \ref{secresnum}. 
This is especially true when using coarse meshes, as is often the case in industrial problems. 
%%%%%%%%%%%% bizarre ?
%%%%%%%%%%% nt in order to obtain a good accuracy on the overall fluxes.  
%%%%%%%%%%%%%%%%%%%%%%%%%
\subsection{\ldots to an optimal compromise?}

Therefore we now propose a scheme which has the advantage of both techniques: we shall use equation \refe{defbilF} and keep the unknowns $u_\edge$ on the edges which require them, for instance those where the matrix $\Lambda$ is discontinuous: hence \refe{consfludis} will hold for all edges associated to these unknowns; for all other interfaces, we shall impose the values of $u$ using \refe{ecrbar}, and therefore eliminate these unknowns.  
Let us decompose the set $\edgesint$ of interfaces into two nonintersecting subsets, that is: 
$\edgesint = \bary\cup \edgeshyb, \edgeshyb=\edgesint \setminus\bary$. 
The interface unknowns associated with $\bary$ will be computed by using the barycentric formula \refe{ecrbar}. 
%%%%%%%%%%% nouveau
\begin{remark}
Note that, although the accuracy of the scheme is increased in practice when the points where the matrix $\Lambda$ is discontinuous are located within the set $\bigcup_{\edge\in \edgeshyb}\edge$, such a property is not needed in the mathematical study of the scheme.
\end{remark}

Let us introduce the space $X_{\disc,\bary} \subset X_{\disc,0}$ defined by:
\begin{equation}
X_{\disc,\bary}=\{v \in X_\disc \mbox{ such that } v_\edge = 0    \mbox{ for all } \edge\in\edgesext   \mbox{ and } v_\edge  \mbox{ satisfying }  \refe{ecrbar}  \mbox{ for all } \edge\in\bary  \}.
\label{defXtilde}
\end{equation}
The composite scheme which we consider in this work reads:
\begin{equation}
\left\{
\begin{array}{l}
\mbox{Find } u \in X_{\disc,\bary} 
%%%% pas clair \mbox{( that is  } (u_K)_{K \in \mesh}, \, (u_\edge)_{\edge \in \edgeshyb})
\mbox{ such that:} \\
\langle u,v \rangle_F   = \sum_{K\in\mesh} v_K \int_K f(\x)\d\x,\ \mbox{ for all } v \in X_{\disc,\bary}.
\end{array}\right.
\label{scheme}
\end{equation}
%%%% il etait temps !!
%%%% provided, of course, that we have an adequate expression of the numerical flux 
%%%%  $F_{K,\sigma}(u) $ with respect to the discrete unknowns. 
We therefore obtain a symmetric scheme with card($\mesh)$ + card$(\edgeshyb)$ equations and unknowns. 
It is thus less expensive while it remains accurate (for the choice of numerical flux given below) even in the case of strong heterogeneity (see section \ref{secresnum}). 

Note that with the present scheme, \refe{consfludis} holds for all $\edge\in\edgeshyb$, but not generally for any $\edge\in\bary$. 
However, fluxes between pairs of control volumes can nevertheless be identified. 
Indeed, we may write
\[
\langle u,v \rangle_F = \sum_{K\in\mesh}\left(\sum_{\edge\in\edgescv\cap \edgeshyb} F_{K,\edge}(u) (v_K - v_\edge)
+ \sum_{\edge\in\edgescv\cap\bary} \sum_{L\in\mesh} F_{K,\edge}(u) \beta_\edge^L (v_K - v_L)\right),
\]
and therefore: 
\[
\langle u,v \rangle_F = \sum_{K\in\mesh}\sum_{\edge\in\edgescv\cap \edgeshyb} F_{K,\edge}(u) (v_K - v_\edge)
+ \half \sum_{ (K,L)\in \mathcal{N}_\disc} F_{K,L}(u) (v_K - v_L),
\]
where 
\[
\mathcal{N}_\disc = \{ (K,L)\in\mesh^2, \exists \edge\in\edgescv\cap\bary, \beta_\edge^L\neq 0 \hbox{ or } 
\exists \edge\in\edgescvv\cap\bary, \beta_\edge^K\neq 0\},
\]
and
\[
F_{K,L}(u) = \sum_{\edge\in\edgescv\cap\bary} F_{K,\edge}(u) \beta_\edge^L -  
\sum_{\edge\in\edgescvv\cap\bary}F_{L,\edge}(u) \beta_\edge^K.
\]
Note that, if $(K,L)\in \mathcal{N}_\disc$, then $(L,K)\in \mathcal{N}_\disc$ and $F_{K,L}(u) = -F_{L,K}(u)$; furthermore, $F_{K,L}(u)\neq 0$ implies $(K,L)\in \mathcal{N}_\disc$, and the scheme's stencil is determined by the set $\{L\in\mesh$ such that  $(K,L)\in \mathcal{N}_\disc\}$.
Then, taking $v_K=1$ and all other degrees of freedom of $v\in X_{\disc,\bary}$ equal to 0, \refe{scheme} yields
\[
\sum_{\edge\in\edgescv\cap \edgeshyb} F_{K,\edge}(u) + \sum_{\substack{L\in\mesh \\ (K,L)\in  \mathcal{N}_\disc}} F_{K,L}(u)
= \int_K f(\x)\d\x,
\]
which shows the ``finite volume philosophy" of the scheme.
\begin{remark}[Other boundary conditions]
In the case of Neumann or Robin boundary conditions, the discrete space $X_{\disc,\bary}$ is modified to include the unknowns associated to the corresponding edges, and the resulting discrete weak formulation is then straightforward. 
\end{remark}
% re nouveau
\begin{remark}[Extension of the scheme]\label{rembarycons}
For consistency reasons, it is preferable that the coefficients $\beta_\sigma^K$ associated with $\edge \in \bary$ be nonzero for points $x_K$ that lie in the same ``regularity zone" of the solution as $x_\edge$ (that is with a zone with no diffusion tensor discontinuity). 
This is not always easy: indeed, in the tilted barrier example described in Section \ref{sec-tilt} below, the barrier contains only one layer of grid cells,  so that, for an internal interface of this layer, it is difficult to use points $x_L$ that are located in the same diffusion regularity zone with respect to $x_K$. 
There is, however, no additional difficulty to replace \refe{ecrbar} in the definition of \refe{defXtilde} by
\begin{equation}
u_\edge = \sum_{K\in\mesh}\beta_\edge^K u_K + \sum_{\edgep\in\edgeshyb}\beta_\edge^\edgep u_\edgep\qquad \forall \edge\in\bary, \label{ecrbarbis}
\end{equation}
\begin{equation}\ba
\dsp\sum_{K\in\mesh}\beta_\edge^K+ \sum_{\edgep\in\edgeshyb}\beta_\edge^\edgep = 1\hbox{ and }
\dsp\xedge = \sum_{K\in\mesh}\beta_\edge^K \xcv+ \sum_{\edgep\in\edgeshyb}\beta_\edge^\edgep\xedgep\qquad \forall \edge \in \bary.
\ea
 \label{propbarybis}
\end{equation}
This trick solves the consistency issue without switching the edge to the hybrid set $\edgeshyb$, while all the mathematical properties shown below still hold.
\end{remark}

\subsection{Construction of the fluxes using a discrete gradient} \label{sec-cons-flux}
For the definition of the schemes to be complete, there now remains to explain how we find a convenient expression for $F_{K,\sigma}(u)$ with respect to the discrete unknowns. 
An idea that has been used in several of the schemes referred to in the Introduction  is to look for a consistent expression of the flux by using adequate linear combinations of the unknowns; however, referring to the beginning of Section \ref{fondam}, such a reconstruction does not in general lead to the desired properties (P2) (symmetric definite positive matrices) and (P3) (convergence). 
Our idea here is different: it is  based on the identification of the numerical fluxes $F_{K,\sigma}(u)$ through
the mesh-dependent bilinear form $\langle \cdot, \cdot \rangle_F$ defined in \refe{defbilF},
using the expression of a discrete gradient.  
Indeed let us assume that, for all $u\in X_{\disc}$, we have constructed a discrete gradient $\nabla_\disc u$, we then seek a family $(F_{K,\edge}(u))_{\stackrel{K \in \mesh}{\edge \in \edgescv}}$ such that 
\begin{equation}
 \langle u,v \rangle_F = \sum_{K\in\mesh} \sum_{\edge\in\edgescv} F_{K,\sigma}(u) (v_K - v_\edge)  = \int_\O \grad_\disc u(\x)\cdot \Lambda(\x)\grad_\disc v(\x)\d\x \qquad \forall u,v\in X_{\disc}.\label{firstidea}
\end{equation}
%re nouveau
\begin{remark}[On the construction of the discrete fluxes] \label{cons-dis-flux}
Note that it is always possible to deduce an expression for $F_{K,\sigma}(u)$ satisfying \refe{firstidea}, under the sufficient condition that, for all $K\in\mesh$ and a.e. $\x\in K$, $\nabla_\disc u(\x)$ is expressed as a linear combination of $(u_\edge - u_K)_{\edge\in\edgescv}$, the coefficients of which are measurable bounded functions of $\x$. 
This property is ensured in the construction of  $\nabla_\disc u(\x)$ given below.
\end{remark}
%re nouveau modifie par rh
We prove in  Section \ref{cvstudy}  below that the desired properties (P2) and (P3) hold if the discrete gradient satisfies the following properties:
\begin{enumerate}
\item (Weak compactness) For a sequence of space discretisations of $\Omega$ with mesh size tending to 0, if the sequence of associated grid functions is bounded in some sense, then their discrete gradient converges at least weakly in $L^2(\O)^d$  to the gradient of an element of $H^1_0(\O)$;
\item (Consistency) If  $\varphi$ is a regular function from $\overline \Omega$ to $\R$, the discrete gradient of the piece-wise function  defined by taking the value $\varphi(\xcv)$ on each control volume $K$ and  $\varphi(\xedge)$ on each edge $\edge$ is a consistent approximation of the gradient of $\varphi$.
\end{enumerate}
Let us first define:
\be
\grad_K u = \frac {1} {\mcv} \sum_{\edge\in\edgescv} \medge (u_\edge - u_K) \ncvedge \qquad \forall K\in\mesh, \forall u\in X_\disc, 
\label{defgradap}\ee
where $\ncvedge$ is the outward to $K$ normal unit vector, $\mcv$ and $\medge$ are the usual measures (volumes, areas, or lengths) of $K$ and $\edge$.
The consistency of formula \refe{defgradap} stems from  the following geometrical relation:
\be
\sum_{\edge\in\edgescv}\medge \ncvedge (\xedge - \xcv)^t = \mcv {\rm Id}  \qquad\forall K\in\mesh,
\label{magical}\ee
where $(\xedge - \xcv)^t$ is the transpose of $\xedge - \xcv \in \R^d$, and $\rm{Id}$ is the $d\times d$ identity matrix.
Indeed, for any linear function defined on $\O$ by  $\psi(\x) = \bG\cdot \x$ with $\bG \in \R^d$, assuming that
 $u_\edge = \psi(\xedge)$ and $u_K = \psi(\xcv)$, we get 
$u_\edge - u_K = (\xedge - \xcv)^t \bG = (\xedge - \xcv)^t\grad \psi$, hence \refe{defgradap} leads to
$\grad_K u = \grad \psi$. 

Since the coefficient of $u_K$ in \refe{defgradap} is in fact equal to zero, a re-construction of the discrete gradient $\nabla_\disc u$ solely based on  \refe{defgradap}  cannot lead to a definite discrete bilinear form in the general case.  
Hence, we now introduce a stabilised gradient: 
\be
\nabla_{K,\edge} u = \nabla_{K} u + R_{K,\edge} u \ \ncvedge,
\label{defgradKedge}\ee
with
\be
R_{K,\edge} u = \frac {\sqrt{d}} {\dcvedge} \left( u_\edge - u_K - \grad_K u\cdot(\xedge-\xcv)\right),
\label{defrkedge}
\ee
(recall that $d$ is the space dimension and $\dcvedge$ is the Euclidean distance between $\xcv$ and $\edge$). 
We may then  define $\nabla_\disc u$ as the piece-wise constant function equal to   $\nabla_{K,\edge} u$ a.e. in the cone $D_{K,\edge}$ with vertex $\xcv$ and basis $\edge$: 
\begin{equation}
\nabla_\disc u (\x) = \nabla_{K,\edge} u  \mbox{ for a.e. } \x \in D_{K,\edge}.
\label{defgrad}
\end{equation}
Note that, from the definition \refe{defrkedge},  thanks to \refe{magical} and to the definition \refe{defgradap}, we get that  
\begin{equation}
 \sum_{\edge \in \edgescv}  \frac {\medge \dcvedge} d  R_{K,\edge} u\ \ncvedge = 0  \qquad\forall K \in \mesh.
\label{propRnulle}
\end{equation}
We prove in Lemmata \ref{consflux} and \ref{consgrad} below that the discrete gradient defined by \refe{defgradap}-\refe{defgrad} indeed satisfies the above stated weak compactness and consistency properties.
In order to identify the numerical fluxes $F_{K,\sigma}(u)$ through Relation \refe{firstidea}, we put the discrete gradient in the form
$$
\nabla_{K,\edge} u = \sum_{\edge' \in \edgescv} (u_{\edge'} - u_\cv) \mathbi{y}^{\edge \edge'},
$$
with 
\be
\mathbi{y}^{\edge  \edge'}= \left\{\begin{array}{ll}
\dsp \frac \medge \mcv \ncvedge +  \frac {\sqrt{d}} {\dcvedge} \left(1 - \frac\medge  \mcv \ncvedge \cdot (\xedge - \xcv)\right)\ncvedge & \mbox{ if } \edge = \edgep \\
\dsp \frac \medgep \mcv \ncvedgep - \frac {\sqrt{d}} {\dcvedge \mcv}\medgep \ncvedgep \cdot (\xedge - \xcv)\ncvedge & \mbox{ otherwise }.
\end{array}\right.
\label{defy}
\ee
Thus,
\be
\int_\O \grad_\disc u(\x)\cdot \Lambda(\x)\grad_\disc v(\x)\d\x = \sum_{K\in\mesh} \sum_{\edge\in\edgescv}\sum_{\edge'\in\edgescv} A_K^{\edge\edge'} (u_\edge - u_K)(v_{\edge'} - v_K)\qquad  \forall u,v\in X_{\disc},
\label{conseq}\ee
with,
\begin{equation}
A_K^{\edge\edge'}  = \sum_{\edge''\in \edgescv} \mathbi{y}^{\edge''  \edge} \cdot \Lambda_{\cv,\edge''}  \mathbi{y}^{\edge''  \edge'} \mbox{ and } \Lambda_{\cv,\edge''} = \int_{D_{K,\edge''}} \Lambda(\x) \d\x.
\label{defA}\end{equation}
Then we get that the local matrices $(A_K^{\edge\edge'})_{\edge\edge'\in\edgescv}$ are symmetric and positive, and the identification of the numerical fluxes using \refe{firstidea} leads to the expression: 
\be 
F_{K,\edge}(u) =  \sum_{\edge'\in\edgescv} A_K^{\edge\edge'} (u_K - u_{\edge'}). \label{defflups}
\ee

\begin{remark}[Link with the MFD method] \label{compmim}
The above technique  yields an explicit construction of  a particular MFD method. 
Indeed, if one chooses $x_K$ as the centre of mass of $K$, the matrix $A_K$ defined by \refe{defA} is an adequate choice for the matrix $\mathbb{W}_E$ which is a parameter in the general formulation of  the family of MFD methods as proposed in \cite{bre-05-fam}.
%ndrh comparer avec le choix de WE que prennent brezzi et al
 The advantages of the specific matrix $A_K$ are that:

\begin{itemize}
\item  on particular meshes, taking a natural choice for $x_K$ (for instance the circumcenter for a triangular mesh in the case of a 2D isotropic problem), it degenerates to a diagonal matrix (see Lemma \ref{superad} below);
\item it is linked to an explicit formulation of a consistent gradient, which is used to define the discrete bilinear form \refe{conseq}.
\end{itemize}
Note however that the  SUSHI scheme defined by \eqref{scheme} is not the MFD method of \cite{bre-05-con,bre-05-fam}; the main reason  is that 
according to the choice $\bary$, the SUSHI scheme may be either a completely cell-centred scheme, or a partly or fully hybrid scheme, while the MFD method is a pure hybrid scheme.
Note also that in SUSHI, one may take any point in cell $K$ for $x_K$, while the MFD schemes \cite{bre-05-con,bre-05-fam} are constructed with the centre of mass (however, this choice might be generalised).

Note that the procedure which we describe in Section \ref{sec-succes} to write a cell-centred scheme could be applied to any mimetic scheme (or low order mixed finite element scheme) to yield a centred scheme.
%and this is also the case for the convergence analysis which is performed in Section \ref{cvstudy}  without any assumption on the regularity of the solution 
However, further investigations are needed to determine under what conditions the present convergence analysis extends to mimetic schemes, and conversely, whether the mimetic analysis applies to the SUSHI scheme  (ongoing work, \cite{dro-09-com}).
\end{remark}

%%%%%% re il faut quand meme dire ce que cette def va permettre
The fluxes defined by \eqref{defgradap}-\eqref{defflups} satisfy certain properties which are detailed in Lemma \ref{casparthyb}, and which allow us to prove the convergence of the scheme, as is shown in  Theorem \ref{cvgce}.
Note that it seems difficult to deduce such properties from fluxes obtained by using natural expansions of regular functions. 
%Then Theorem \ref{cvgce} shows that these properties are sufficient to provide the convergence of the scheme. 
Note also that  both Lemma  \ref{casparthyb} and Theorem \ref{cvgce} hold for general heterogeneous, anisotropic and possibly discontinuous fields $\Lambda$, for which the solution $u$ of \refe{ellgenf} is not in general more regular than $u\in H^1_0(\O)$. 
In the case where $\Lambda$ and $u$ are regular enough, the local  flux consistency satisfied by \refe{defflups} is used in order to obtain an error estimate, see Theorem \ref{eresthyp}.
%%%
The coefficient $\sqrt d$ may be replaced by any positive real number without any change in the proof of convergence; in fact, for certain problems it can be interesting to use another coefficient, as described in \cite{eym-08-ben} for the so called ``SUSHI-P" scheme  (P for parametric, meaning that the user may choose the stabilisation coefficient as well as the set of edges $\bary$). 
The choice $\sqrt d$ is however natural in the sense that with this value, if $\bary = \emptyset$,  the scheme boils down in two dimensions to the well-known harmonic averaging five points scheme on rectangles and a four-point scheme on triangles; more generally, in any space dimension, even if $\bary \not = \emptyset$ and taking the most natural value for $u_\sigma$ if $\sigma \in \bary$, the resulting flux is a two-point flux on meshes that satisfy the ``superadmissibility condition" \eqref{propsuperad}, not necessarily with a harmonic averaging of $\Lambda$ in the case $\edge \in \bary$; this is  proven in the next lemma. 
Note that this superadmissibility condition is also satisfied by rectangular parallelepipeds in three dimensions but unfortunately not by tetrahedra.
\begin{lemma}[Superadmissible mesh and two-point flux]\label{superad}
Let $\disc=(\mesh,\edges,\points)$ be a discretisation of $\Omega$ in the sense of Definition \ref{adisc},  satisfying the following superadmissibility  condition:
\begin{equation}
\ncvedge = \frac {\xedge - \xcv} {\dcvedge} \qquad\forall K \in \mesh, \; \forall \edge \in \edgescv.
\label{propsuperad}
\end{equation} 
Let us furthermore assume that $\Lambda(\x) = \lambda(\x){\rm Id}$,  where $\lambda$ is a piece-wise constant function from $\Omega$ to $\R$, which is equal to a constant $\lambda_K$ in each $K\in\mesh$; then, the inner product defined by \eqref{firstidea}-\eqref{defrkedge} reads:
$$ 
\left\langle u,v\right\rangle_F = \sum_{K\in\mesh}  \lambda_K \sum_{\edge\in\edgescv} \frac \medge \dcvedge (u_K - u_\edge) (v_K - v_\edge).
$$
Moreover, choosing thanks to \refe{propsuperad}, $\xedge = (\dcvedge \xcvv + \dcvvedge \xcv)/(\dcvedge+ \dcvvedge)$ for $\edge\in\edgesint$ with
$\mesh_\edge =\{K, L\}$ in  \refe{propbary}, the scheme \eqref{scheme} is the following two-point flux scheme:
\begin{eqnarray}
\sum_{K \in \mesh} F_{K,\sigma} &=& \int_K f(\x) \ \d \x,\\
F_{K,\sigma} & = & 
\frac {  {\lambda_K} {\lambda_L} (\dcvedge + \dcvvedge)} 
{  {\lambda_K} {\dcvvedge} +   {\lambda_L} {\dcvedge}} \frac {\medge} {\dcvedge + \dcvvedge}( u_K - u_L ) \mbox{ if } \edge \in\edgesint\cap \edgeshyb,\ \mesh_\edge =\{K, L\},\\
F_{K,\sigma} & =&\frac {\dcvedge \lambda_K + \dcvvedge \lambda_L} {\dcvedge + \dcvvedge}
\frac {\medge} {\dcvedge + \dcvvedge}( u_K - u_L) \mbox{ if } \edge \in\edgesint\cap \bary,\ \mesh_\edge =\{K, L\},\\
F_{K,\sigma} & =& \lambda_K \frac \medge \dcvedge u_K \mbox{ if } \edge  \in \edgesext \cap \edgescv,\end{eqnarray}
\end{lemma}
\begin{proof}
Let us compute $\left\langle  u,v\right\rangle_F$ under the assumptions of Lemma \ref{superad}. 
From \refe{firstidea} and thanks to  \refe{propRnulle} we get: 
\begin{eqnarray}
\left\langle  u,v\right\rangle _\disc & = & \sum_{K\in\mesh} \lambda_K  \int_K \grad_\disc u(\x)\cdot \grad_\disc v(\x)\d\x \nonumber\\
 & = &  \sum_{K\in\mesh} \lambda_K \left(  \mcv \nabla_{K} u\cdot \nabla_{K} v
+  \sum_{\edge\in\edgescv} \frac {\medge\dcvedge} {d} \ R_{K,\edge} u \ R_{K,\edge} v\right). \nonumber
\end{eqnarray}
Now from the definition \eqref{defgradap} and thanks to the assumption \eqref{propsuperad}, the discrete gradient given by \refe{defgradap} may be written as follows:
$$
\grad_K v = \frac {1} {\mcv} \sum_{\edge\in\edgescv} \frac \medge  \dcvedge (v_\edge - v_K) (\x_\edge - \x_K)\qquad \forall K\in\mesh,\ \forall v\in X_\disc, 
$$
From \refe{magical}, we get
\[
\sum_{\edge\in\edgescv} \frac {\medge\dcvedge} {d} \frac {\sqrt{d}} {\dcvedge} 
(\xedge-\xcv) \ \frac {\sqrt{d}} {\dcvedge} (\xedge-\xcv)^t  = \mcv {\rm Id}.
\]
Therefore, we get that
\[
\dsp\sum_{\edge\in\edgescv} \frac {\medge\dcvedge} {d} \ R_{K,\edge} u \ R_{K,\edge} v = \dsp
\sum_{\edge\in\edgescv} \frac {\medge} {\dcvedge} ( u_\edge - u_K)( v_\edge - v_K)  - \mcv \nabla_{K} u\cdot \nabla_{K} v,
\]
which in turn yields that
\[
%\int_\O \grad_\disc u(\x)\cdot \Lambda(\x)\grad_\disc v(\x)\d\x 
\left\langle  u,v\right\rangle _\disc= \sum_{K\in\mesh} \lambda_K 
\sum_{\edge\in\edgescv} \frac {\medge} {\dcvedge} ( u_\edge - u_K)( v_\edge - v_K).
\]
Hence the matrix  $A_K$ only contains the terms  $\frac {\medge} {\dcvedge}$ on the diagonal, and the flux $F_{K,\edge}(u)$ is given by
\[
F_{K,\edge}(u) =  \lambda_K \frac {\medge} {\dcvedge} ( u_K - u_\edge).
\]
Then the scheme \refe{scheme} can be written as a classical cell-centred  finite volume scheme, with two-point fluxes
$F_{K,L}(u) = -F_{K,L}(u)$ for any $\edge \in\edgesint$ with $\mesh_\edge = \{K,L\}$.
Indeed, in the case $\edge\notin\bary$,  
the above expression of $F_{K,\edge}(u)$ allows us to get the following expression of $u_\edge$ from \refe{consfludis}:
\[
u_\edge = \frac {\frac {\lambda_K} {\dcvedge} u_K + \frac {\lambda_L} {\dcvvedge} u_L}{\frac {\lambda_K} {\dcvedge} + \frac {\lambda_L} {\dcvvedge}}.
\]
This yields the harmonic averaging two-point flux
\[
F_{K,L}(u) = 
\medge\frac {\frac {\lambda_K} {\dcvedge} \frac {\lambda_L} {\dcvvedge}} 
{\frac {\lambda_K} {\dcvedge} + \frac {\lambda_L} {\dcvvedge}} ( u_K - u_L ).
\]
In the case $\edge\in\bary$, the two-point barycentric formula $u_\edge = (\dcvedge u_\cvv + \dcvvedge u_\cv)/(\dcvedge+ \dcvvedge)$ together with 
 \refe{scheme}  leads to the resulting  two-point flux  
\[
F_{K,L}(u) = \frac {\dcvedge \lambda_K + \dcvvedge \lambda_L} {\dcvedge + \dcvvedge}
\frac {\medge} {\dcvedge + \dcvvedge}( u_K - u_L).
\]
\end{proof}

\section{Numerical results}\label{secresnum}

We present some numerical results obtained with various choices of $\bary$ in the scheme \refe{scheme}, \refe{defbilF} with the flux \refe{defflups}, which we synthesise here for the sake of clarity: 
 \be 
 \left\{\begin{array}{l}
\dsp \mbox{Find } u \in X_{\disc,\bary} \mbox{ (that is  } (u_K)_{K \in \mesh}, \, (u_\edge)_{\edge \in \edgeshyb}),  \mbox{ such that:} \\
\dsp \sum_{K\in\mesh} \sum_{\edge\in\edgescv} F_{K,\sigma}(u) (v_K - v_\edge)  = \sum_{K\in\mesh} v_K \int_K f(\x)\d\x,\ \mbox{ for all } v \in X_{\disc,\bary}, \\
\dsp \mbox{with } F_{K,\edge}(u) =  \sum_{\edge'\in\edgescv} A_K^{\edge\edge'} (u_{\edge'} - u_K)\qquad \forall K \in  \mesh, \forall \edge \in \edgescv.
\end{array}\right.
 \label{schemepratique}
 \ee
where the matrices  $A_K^{\edge\edge'}$ are defined by \eqref{defA}-\eqref{defy}.
In the following, we shall use the choices $\bary = \emptyset$ (HFV), $\bary = \edgesint$ or $\bary$ the set of edges which are located on the diffusion tensor discontinuity interfaces; this latter choice is reported as SUSHI-NP (for non parametric) in \cite{eym-08-ben}, in contrast with SUSHI-P (for parametric) where the choice of the set $\bary$ may be different, along with the value of the stabilisation coefficient in \eqref{defrkedge}.

\subsection{Implementation}\label{numimp}

Let us first describe an implementation aspects of the scheme.
The unknowns, {\it i.e.} the values $u_K$, for $K\in\mesh$ and the values $u_\edge$, $\edge\in\edgesint\cap \edgeshyb$, are ordered as $(u_i)_{i=1,\dots,N}$.
The $N\times N$ matrix and the $N\times 1$ right-hand-side of the  linear system resulting from \eqref{scheme} are computed thanks to a loop over the control volumes $K\in\mesh$ and to an inner loop on each edge $\edge\in\edgescv$.
Let us detail the matrix computation loop.
\begin{enumerate}
\item All stored matrix coefficients are initially set to 0.
\item The expression $F_{K,\edge}(u)$ is written in the form  $F_{K,\edge}(u) = \sum_{i=1,\ldots,N} a_{K,\edge}^{(i)} u_i$, where the nonzero coefficients
$(a_{K,\edge}^{(i)})_{i=1,\ldots,N}$ are only locally computed (they are not stored for all $K$ and $\edge$). These coefficients are obtained after the elimination  of all $(u_\edge)_{\edge\in\edgescv\cap\bary} $ in \refe{defflups}: 
\[
F_{K,\edge}(u) =  \sum_{\edge'\in\edgescv\cap \edgeshyb} A_K^{\edge\edge'} (u_K - u_{\edge'}) +
\sum_{\edge'\in\edgescv\cap\bary} A_K^{\edge\edge'} \sum_{L\in\mesh} \beta_{\edge'}^L  (u_K - u_L).
\]
\item The line of the matrix corresponding to the unknown $u_K$ is incremented at the column $j$ with the coefficient $a_{K,\edge}^{(j)}$.
\item If $\edge\in \bary$ with $v_\edge = \sum_{L\in\mesh}\beta_\edge^L v_L$ for any $v\in X_{\disc,\bary}$, the line of the matrix corresponding
to each $L\in\mesh$ such that $\beta_\edge^L\neq 0$ is incremented at the column $j$ with the coefficient
$- \beta_\edge^L a_{K,\edge}^{(j)}$.
\item If  $\edge\in \edgesint\cap \edgeshyb$, the line of the matrix corresponding to the edge $\edge$   is incremented
at the column $j$ with the coefficient $-a_{K,\edge}^{(j)}$.
\end{enumerate}

This procedure is identical in the cases $\bary = \emptyset$ (HFV),  $\bary \not = \emptyset$  and $\bary = \edgesint$. 
However, in the case where $\bary=\emptyset$ (HFV), one may eliminate the unknowns $u_K$ with respect to the unknowns $u_\edge$, as in the hybrid implementation of the mixed finite element method.

\subsection{Order of convergence}
We consider here the numerical resolution of Equation \refe{eq1} supplemented by the homogeneous Dirichlet boundary condition \eqref{cl0};  the right-hand side is chosen so as to obtain an exact solution to the problem and easily compute the error between the exact and approximate solutions. 
We consider Problem  \refe{eq1}-\refe{cl0} with a constant matrix $\Lambda$:
\be 
\Lambda = \left( \begin{array}{cc}
1.5 & .5 \\ .5 & 1.5
\end{array}
\right)  \label{clubar}, 
\ee
and choose $f:$  $\Omega \to \R$ such that the exact solution to Problem  \refe{eq1}--\eqref{cl0} is $\bar u$ defined by $\bar u(x,y) = 16 x(1-x) y(1-y)$ for any $(x,y) \in \overline{ \Omega}$.  
Note that in this case, the composite scheme is in fact the cell-centred scheme, there are no edge-unknowns.  

Let us first consider conforming meshes, such as the triangular meshes which are depicted on Figure \ref{fig_mesh1}, and uniform square meshes.  
\twofig{mesh1_1}{mesh1_4}{Regular conforming coarse and fine triangular grids}{fig_mesh1}

For both $\bary = \emptyset$ (pure hybrid scheme: HFV) and $\bary = \edgesint$ (cell-centred scheme), the order of convergence is close to 2 for the unknown $u$ and 1 for its gradient. 
Of course, the hybrid scheme is almost three times more costly in terms of number of unknowns than the cell-centred scheme  for a given precision.  
However, the number of nonzero terms in the matrix is, again for a given precision on the approximate solution, larger for the cell-centred scheme than for the hybrid scheme. 
Hence the number of unknowns is probably not a sufficient criterion for assessing the cost of the scheme.

Results were  also obtained in the case of uniform square or rectangular meshes. They show a better rate of convergence of the gradient (order 2 in the case of $\edgeshyb =\edgesint$  and 1.5 in the case  $\bary =\edgesint$), even though the rate of convergence of the approximate solution remains unchanged and close to 2.

We then use a rectangular nonconforming mesh, obtained by cutting vertically the domain into two parts and using a rectangular grid of $3n \times 2n$ (resp. $5n \times 2n$) on the first (resp. second side), where $n$ is the number of the mesh, $n = 1, \ldots, 7$.  
Again, the order of convergence which we obtain is $2$ for $u$ and around $1.8$ for the gradient. 
We give in Table \ref{table_cnc} below the errors obtained in the discrete $L^2$ norm for $u$ and $\nabla u$ for a nonconforming mesh and  (in terms of number of unknowns) and for the rectangular $4\times 6$  and $4 \times 10$ conforming rectangular meshes, for both the hybrid and cell-centred schemes. 
We show in Figure \ref{figprofilpression} the solutions for the corresponding grids (which look much the same for the two schemes). 

\begin{center}
\begin{table}

\begin{center}
\begin{tabular}{|c|cc|cc|cc|cc|cc|}
  \hline
   & \multicolumn{2}{c|}{NU} & \multicolumn{2}{c|}{NM} & \multicolumn{2}{c|}{$\epsilon(u)$}   &\multicolumn{2}{c|}{$\epsilon(\nabla u)$}   \\
 \hline
 n   &  $\bary =\emptyset$ & $\bary =\edgesint$ &  $\bary =\emptyset  $ &  $\bary =\edgesint $ &   $\bary =\emptyset  $  &    $\bary =\edgesint $&      $\bary =\emptyset  $  &  $\bary =\edgesint $ \\
\hline
C1&     130&  48  & 874 & 488  & 1.28E-01&  1.20E-01&  1.64E-02&  3.57E-02\\
NC&     182&  64  & 1334& 724  & 1.03E-01&  9.43E-02&  1.66E-02&  3.69E-02\\
C2&     222&  80  & 1542& 864  & 7.61E-02&  7.09E-02&  9.18E-03&  2.44E-02\\
% 1&     130&     874&   7.11E-15&   1.28E-01&   1.64E-02\\
% 2&     182&    1334&  5.33E-15&  1.03E-01&  1.66E-02\\
% 3&     222&    1542&  2.04E-14&  7.61E-02&  9.18E-03 \\
\hline
\end{tabular}
\end{center}
\caption{\small Error for the nonconforming rectangular mesh, pure hybrid scheme ($\bary = \emptyset$) and centred ($\bary = \edgesint$) schemes. 
For both schemes NU is the number of unknowns in the resulting linear system, NM is the number of nonzero terms in the matrix, $\epsilon(u)$ is the discrete $L^2$ norm of the error of the solution and $\epsilon(\nabla u)$ is the discrete $L^2$ norm of the error in the gradient. C1 and C2 are the two conforming meshes represented on the left and the right in Figure \ref{figprofilpression}, and NC is the nonconforming one represented in the middle.}
\label{table_cnc}

\end{table}
\end{center}

\threefig{hybrid_6_6}{hybrid_6_10}{hybrid_10_10}{The approximate solution for conforming and nonconforming meshes. Left: conforming $8\times6$ mesh, centre: nonconforming $4\times 6,4\times 10$ mesh, right: conforming $10\times10$. }{figprofilpression}

Further detailed results on several problems and conforming, nonconforming and distorted meshes may be found in \cite{eym-08-ben}.

\subsection{The case of a highly heterogeneous tilted barrier} \label{sec-tilt}
\begin{figure}[ht]\begin{center}
    \scalebox{.52}{\input{dom.\sufftex}}\hskip 0.4cm \includegraphics[width=3.5cm]{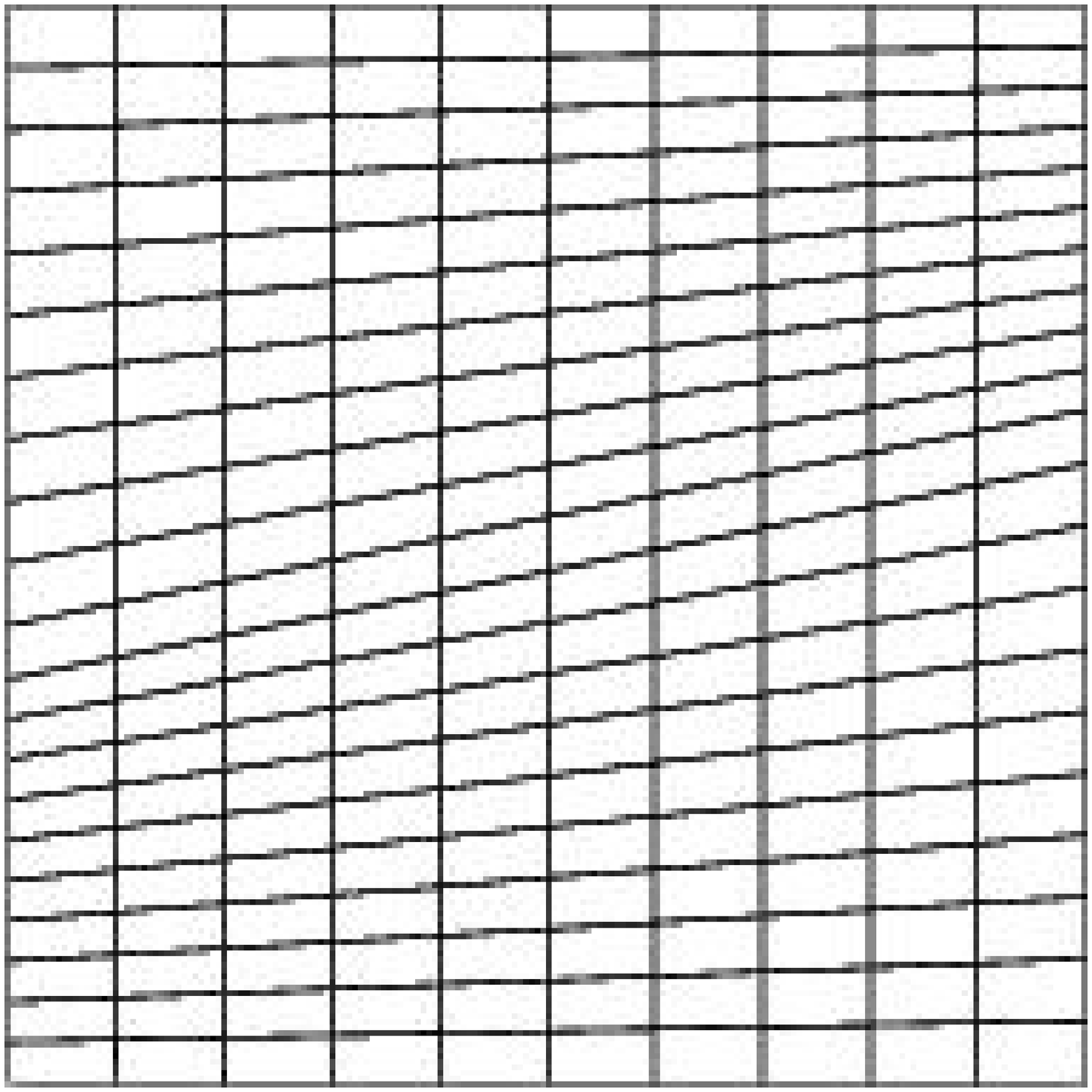}\hskip 0.4cm \includegraphics[width=3.5cm]{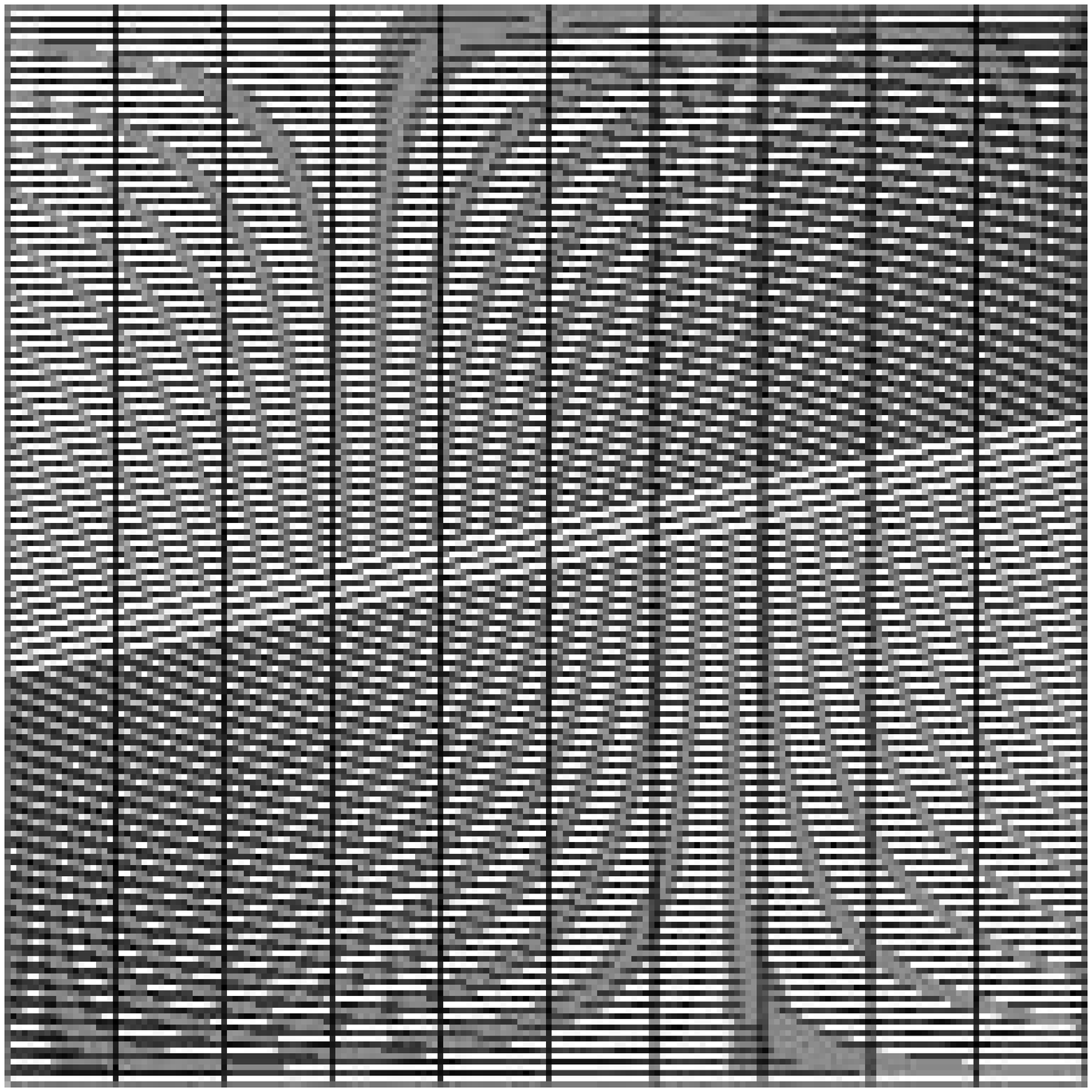}
    \end{center}
\caption{Domain and meshes used for the tilted barrier test: mesh 1 ($10\times 21$ centre), mesh 2 ($10\times 100$ right)\label{fig.mesh5}}\end{figure}

We now turn to the heterogeneous case. 
The domain $\Omega = ]0,1[\times]0,1[$ is composed of 3 sub-domains, which are depicted in Figure \ref{fig.mesh5}:  
$\Omega_1 = \{(x,y) \in \Omega; \phi_1(x,y) <0 $\}, with $\phi_1(x,y)= y - \delta (x - .5) - .475$,
$\Omega_2 = \{(x,y) \in \Omega; \phi_1(x,y) >0, \phi_2(x,y) <0 $\}, with $\phi_2(x,y)= \phi_1(x,y) - 0.05$,     $\Omega_3 = \{(x,y) \in \Omega; \phi_2(x,y) >0$\},
and $\delta = 0.2$ is the slope of the drain (see Figure \ref{fig.mesh5}). 
Dirichlet boundary conditions  are imposed by setting  the boundary values to those of the analytical solution given by $u(x,y) =  -\phi_1(x,y)$ on $\Omega_1\cup \Omega_3$ and  $u(x,y) =  -\phi_1(x,y)/10^{-2} $ on $\Omega_2$.

The permeability tensor $\Lambda$ is heterogeneous and isotropic, given by $\Lambda(\x) = \lambda(\x){\rm Id}$, 
with $\lambda(\x)  = 1$ for a.e. $x\in\Omega_1\cup \Omega_3$ and  $\lambda(\x)  = 10^{-2}$ for a.e. $x\in\Omega_2$.
Note that the isolines of the exact solution are parallel to the boundaries of the sub-domain, and that the tangential component of the gradient is 0.
We use the meshes depicted in Figure \ref{fig.mesh5}. Mesh 3 (containing $10\times 25$ control volumes)
is obtained from Mesh 1 by the addition of two layers of very thin control volumes around each of the two lines of discontinuity of $\Lambda$: because of the very low thickness of these layers, equal to $1/10000$, the picture representing Mesh 3 is not different from that of Mesh 1.

%\threefig{dom}{grid}{gridfin}{Domain and meshes used for the tilted barrier test: mesh 1 ($10\times 21$ center), mesh 2 ($10\times 100$ right)}{fig.mesh5}
%%%%% RE je ne comprends pas 10x100 ??????????
%%%%%%%%%% je ne comprends pas non plus le 1000 unknowns de succes mesh 2

We get the following results for the approximations of the four fluxes at the boundary.

\begin{center}
\begin{tabular}{|c c|c|c|c|c|c|c|}
\hline
  &  & nb. unknowns & matrix size &$x=0$ & $x=1$ & $y=0$ & $y=1$ \\
\hline 
analytical& &  & & -0.2 &  0.2 & 1. & -1.\\
\hline 
\multirow{3}{*}{$\bary =\edgesint$}
%{\begin{sideways}{centered}\end{sideways}} 
% \multirow{3}{*}{\rotatebox{90}{ \mbox{centered}} } 
&mesh 1 & 210 & 2424 & $-1.17$ &  1.17 & 3.51 & $-3.51$\\%OK
%\hline 
&mesh 2 & 1000 & 11904 &$ -0.237$ &  0.237 &  1.104 & $-1.104$ \\ %OK
%\hline 
&mesh 3 & 250 & 2904 & $-0.208 $&  0.208 & 1.02 &$ -1.02$\\
\hline 
\multirow{2}{*}{SUSHI-NP}
 &mesh 1& 239 & 2583 & $-0.2$ &  0.2 & 1. & $-1.$\\
 & mesh 2 &1020&   12036  & $-0.2$ &  0.2 & 1. &$ -1.$\\
\hline 
\multirow{2}{*}{HFV}
 & mesh 1 & 599 & 4311 & $-0.2$ &  0.2 & 1. & $-1.$\\
  & mesh 2 &  2890&   21138& $ -0.2 $&  0.2 & 1. & $-1.$\\
\hline 
\end{tabular}
\end{center}
Note that the values of the numerical solution given by the pure hybrid (HFV) and composite (SUSHI-NP) schemes are equal to those of the analytical solution (this holds under the only condition that the interfaces located  on the lines $\phi_i(x,y)=0$, $i=1,2$, are not included in $\bary$, and that, for all $\edge\in\bary$, 
all $K\in\mesh$ with $\beta_\edge^K\neq 0$ are included in the same  sub-domain $\O_i$). 
Note that Mesh 3, which leads to acceptable results for the computation of the fluxes, is not well suited for such a coupled problem, because of too small control volume measures. 
Hence SUSHI on Mesh 1 appears to be the most suitable method for this problem.

A satisfying natural choice (SUSHI-NP in the above results) is thus to match $ \edgeshyb$ with the discontinuities of $\Lambda$. 
It is sometimes interesting to choose another set $\bary$. 
This is for instance the case for the numerical locking problem for which the choice $\bary = \emptyset$ is best even though the diffusion tensor is homogeneous \cite{eym-08-ben}.

It is also sometimes interesting to replace the stabilisation coefficient $ \sqrt d$ in \eqref{defrkedge} by some other coefficient $\alpha>0$. 
This is the case for instance for very distorted meshes or singular problems, in order to maintain the positivity of the unknown. 
The coefficient $\alpha$ is taken to be greater than $\sqrt d$. 
The approximate solution remains positive, but the $L^2$ norm of the error is generally larger.
We refer to \cite{eym-08-ben} for such experiments.

\section{Convergence of the scheme}\label{cvstudy}

Let us first introduce some notations related to the mesh. 
Let $\disc=(\mesh,\edges,\centers)$ be a discretisation of $\O$ in the sense of Definition \ref{adisc}.
The size of the discretisation $\disc$ is defined by:
\[
h_\disc= \sup\{h_K, K\in \mesh\},
\]
and  the regularity of the mesh  by:
\be
\theta_\disc = \max\left(\max_{\edge\in\edgesint, K,L \in  \mesh_\edge} \frac {\dcvedge} {\dcvvedge}, \max_{K \in \mesh, \edge \in \edgescv} \dfrac{h_K}{\dcvedge} \right).
\label{defetadisc}
\ee

For a given set  $\bary \subset \edgesint$ and for a given family $(\beta_\edge^K)_{\stackrel{K\in\mesh}{\edge\in\edgesint}}$  satisfying property \refe{propbary}, we introduce a measure of the resulting regularity by
\begin{eqnarray}
\theta_{\disc,\bary} =\max\left(\theta_\disc,\max_{K\in\mesh,\edge\in\edgescv\cap\bary} 
\frac  {\sum_{L\in\mesh} \vert\beta_\edge^L\vert |\xcvv - \xedge|^2 } {h_K^2}\right).
\label{regutousazimut}\end{eqnarray}
%%%%%%%%%%%%% nouveau
\begin{remark}
Note that, for any mesh, it is easy to choose the family $(\beta_\edge^K)_{\stackrel{K\in\mesh}{\edge\in\edgesint}}$  so that $\theta_{\disc,\bary}$ remains small. 
It suffices to express $\xedge$ as the barycentre of $d+1$ points $\xcvv$ (which is always possible), for $\cvv$ sufficiently close to $K$, so that $\xcvv - \xedge$ is close to $h_K$ when $\beta_\edge^K \ne 0$. 
% re je ne vois pas ce que tu as voulu dire ici 
% In practice, in the 2D case, we use the three points $\xcvv$ which are closest to $\xedge$.
Note also that in fact, it would be sufficient to have $h_K^\eta$ with $\eta > 1$ instead of $h_K^2$ in  \refe{regutousazimut} thus allowing the use of farther points.
\end{remark}
Remark that, thanks to the assumption that  $K$ is $\xcv$-star-shaped, the following property holds:
\be
\sum_{\edge\in\edgescv} \medge d_{K,\edge} = d\ \mcv\qquad \forall K\in\mesh.
\label{summdiam}\ee

The space $ X_{\disc}$ defined in \refe{defX} is equipped with the following semi-norm: 
\be
\forall v\in X_{\disc}, \ \vert v\vert_{X}^2 = \dsp \sum_{K\in\mesh}  \sum_{\sigma\in\edgescv}
\frac {\medge}{\dcvedge} (v_\edge - v_K)^2,
\label{amudis}\ee
which is a norm on the spaces $X_{\disc,0}$ and  $X_{\disc,\bary}$ respectively defined by \refe{defX0} and \refe{defXtilde}.

Let $H_\mesh(\O)\subset L^2(\O)$  be the set of piece-wise constant functions on the control volumes of the mesh $\mesh$. 
We then denote, for all $v\in H_\mesh(\O)$ and for all $\edge\in\edgesint$ with $\mesh_\sigma= \{K,L\}$,  
$D_\sigma v=\vert v_K-v_L \vert$ and $d_\sigma = \dcvedge + \dcvvedge$, 
and for all $\edge\in\edgesext$ with  $\mesh_\sigma= \{K\}$, we denote  $D_\sigma v=\vert v_K \vert$ 
and $d_\sigma = \dcvedge$.
We then define the following norm:
\be
\forall v\in  H_\mesh(\O),\ \Vert v \Vert_{1,2,\mesh}  =\sum_{K \in \mesh} \sum_{\sigma \in \edgescv} \medge d_{K,\sigma} \left(\frac {D_\sigma v} {d_\edge}\right)^2 = \sum_{\sigma \in \edges} \medge \frac {(D_\sigma v)^2} {d_\edge}.
\label{normtiero}\ee
(Note that this norm is also defined  by \refe{defnorunpd} in Lemma \ref{sobp}, setting $p=2$).

For all $v\in X_{\disc}$, we denote by $\Pi_\mesh v\in H_\mesh(\O)$ the piece-wise function from $\Omega$ to $\R$ defined by $\Pi_\mesh v(\x) = v_K$ for a.e. $\x\in K$, for all $K\in\mesh$. 
Using the Cauchy-Schwarz inequality, we have for all $\edge\in\edgesint$ with $\mesh_\sigma= \{K,L\}$,  
$$\frac {(v_K - v_L)^2} {d_\edge} \le  \frac {(v_K - v_\edge)^2} \dcvedge + \frac {(v_\edge - v_L)^2} \dcvvedge \qquad \forall v\in X_{\disc},$$ 
which leads to the relation
\be
\Vert \Pi_\mesh v \Vert_{1,2,\mesh}^2 \le    \vert v\vert_{X}^2\qquad \forall v\in X_{\disc,0}.
\label{comparnorm}\ee

 For all $\varphi \in C(\Omega, \R)$, we denote by $P_\disc \varphi$ the element of  $X_{\disc}$ 
defined by $((\varphi(\xcv))_{K\in\mesh},(\varphi(\xedge))_{\edge\in\edges})$, by 
 $P_{\disc,\bary} \varphi$ the element $v\in X_{\disc,\bary}$ such that $v_K = \varphi(\xcv)$ for all $K\in\mesh$,
$v_\edge = 0$ for all $\edge\in\edgesext$, $v_\edge = \sum_{K\in\mesh} \beta_\edge^K \varphi(\xcv)$
for all $\edge\in\bary$ and $v_\edge = \varphi(\xedge)$ for all $\edge\in\edgeshyb$.

We denote by $P_\mesh \varphi\in H_\mesh(\O)$ the function such that
$P_\mesh \varphi(\x) = \varphi(\xcv)$ for a.e. $\x\in K$, for all $K\in\mesh$ (we then have
$P_\mesh \varphi = \Pi_\mesh P_\disc\varphi = \Pi_\mesh P_{\disc,\bary} \varphi$).

The following lemma provides an equivalence property between the $L^2$-norm of the discrete gradient, defined by 
\refe{defgradap}-\refe{defgrad} and the norm $\vert\cdot\vert_X$.

\begin{lemma}\label{lemgrad}
Let $\disc$ be a discretisation of $\Omega$ in the sense of Definition \ref{adisc}, and let $\theta \ge \theta_\disc$
be given (where $\theta_\disc$ is defined by \refe{defetadisc}). 
Then there exists $\ctel{101}>0$ and $\ctel{102}>0$ only depending on $\theta$ and $d$ such that:
\begin{equation}
\cter{101} \vert u \vert_X \le \Vert \nabla_\disc u \Vert_{L^2(\Omega)} \le \cter{102} \vert u \vert_X \qquad
\forall u\in X_\disc,
\label{ineqgrad}\end{equation}
where $\nabla_\disc$ is defined by \refe{defgradap}-\refe{defgrad}.
\end{lemma}
\begin{proof}
By definition, 
$$
\Vert \nabla_\disc u \Vert_{L^2(\Omega)^d}^2 = \sum_{K \in \mesh} \sum_{\edge \in \edgescv} \frac {\medge \dcvedge} d \vert \nabla_{K,\edge} u \vert^2. 
$$
Therefore, using property \refe{propRnulle},
\begin{equation}
\Vert \nabla_\disc u \Vert_{L^2(\Omega)^d}^2 = \sum_{K \in \mesh} \left(\mcv \vert \nabla_K u \vert^2 + \sum_{\edge \in \edgescv} \frac {\medge \dcvedge} d   (R_{K,\edge} u)^2\right).
\label{gradR}
\end{equation}

Let us now notice that the following inequality holds:
\be
(a-b)^2 \ge \frac {\lambda} {1 + \lambda} a^2 - \lambda b^2\qquad \forall a,b\in\R,\ \forall \lambda >-1.
\label{coercineq}\ee
%(it is an immediate consequence of $\left(a\sqrt{\frac {1} {1 + \lambda}}  - b\sqrt{1 + \lambda}\right)^2 \ge 0$).
We apply this inequality  to $\left(R_{K,\edge} u \right)^2$ for some  $\lambda > 0$ and obtain
\be
\left( R_{K,\edge} u \right)^2 \ge \frac {\lambda d} {1 + \lambda} \left( \frac {u_\edge - u_K}{\dcvedge} \right)^2 - \lambda d |\grad_K u|^2 \left( \frac {|\xedge-\xcv|}{\dcvedge} \right)^2.
\label{inega}\ee
This leads to
\[
\sum_{\edge\in\edgescv} \frac {\medge\dcvedge} d \left( R_{K,\edge} u \right)^2 \ge \frac {\lambda} {1 + \lambda} \sum_{\edge\in\edgescv}  {\medge\dcvedge} \left(\frac {u_\edge - u_K}{\dcvedge} \right)^2 - 
\lambda \ \mcv  \ d |\grad_K u|^2 \theta^2.
\]
Choosing   $ \lambda =\frac  {{1}} {d \theta^2}, $ we get that
\[
\Vert \nabla_\disc u \Vert_{(L^2(\Omega))^d}^2 \ge \frac {\lambda} {1 + \lambda}  \vert u \vert_X^2,
\]
which shows the left inequality of \refe{ineqgrad}.

Let us now prove the right inequality. On one hand, using the definition \refe{defgradap} of $\nabla_K u$ and \refe{summdiam}, 
the Cauchy--Schwarz inequality leads to

\begin{equation}
\vert \nabla_K u \vert^2 \le \ \frac 1 { {\mcv}^2}\sum_{\edge \in \edgescv} \frac \medge \dcvedge (u_\edge - u_\cv)^2 \sum_{\edge \in \edgescv}  \medge \dcvedge = \frac d \mcv \sum_{\edge \in \edgescv}  \frac \medge \dcvedge (u_\edge - u_\cv)^2. 
\label{etun}
\end{equation}

On the other hand, by definition \refe{defrkedge}, and thanks to the definition of the regularity of the mesh \refe{defetadisc}, we have
\begin{equation}
(R_{\cv,\edge} u)^2 \le 2 d \left( (\dfrac{u_\edge - u_K}{\dcvedge})^2 + \vert \nabla_K u \vert^2 \vert\dfrac{ \xedge - \xcv}{\dcvedge}\vert^2 \right) \le 2 d  \left( (\dfrac{u_\edge - u_K}{\dcvedge})^2+ \theta^2 \vert \nabla_K u \vert^2  \right).
\label{etdeux}
\end{equation}

{F}rom \refe{gradR}, \refe{etun} and \refe{etdeux}, we conclude that the right inequality of \refe{ineqgrad} holds. 
\end{proof}

We may now state a  weak compactness result for the discrete gradient.
\begin{lemma}[Weak discrete $H^1$ compactness]\label{consflux}
Let ${\mathcal F}$ be a family of discretisations in the sense of Definition \ref{adisc} such that there exists $\theta>0$ with $\theta \ge \theta_\disc$ for all $\disc \in {\mathcal F}$.   
Let $(u_\disc)_{\disc\in\mathcal F}$  be a family of functions, such that:
\begin{itemize}
\item $u_\disc\in X_{\disc,0}$ for all $\disc\in\mathcal F$,
\item there exists $C>0$ with $\vert u_\disc\vert_X \le C$ for all $\disc\in\mathcal F$,
\item there exists $u\in L^2(\O)$ with $\lim\limits_{h_\disc \to 0}\Vert \Pi_\mesh u_\disc - u\Vert_{L^2(\O)} = 0$. \end{itemize}
Then, $u\in H^1_0(\O)$ and $\grad_\disc u_\disc$ weakly converge in $L^2(\O)^d$ to $\grad u$ as $h_\disc\to 0$, where the operator $\nabla_\disc$ is defined by \refe{defgradap}-\refe{defgrad}. 
\end{lemma}
\begin{proof}
Let us prolong $\Pi_\mesh u_\disc$ and $\grad_\disc u_\disc$ by 0 outside of $\O$.
Thanks to Lemma \ref{lemgrad},  up to a subsequence, there exists some function $\bG\in L^2(\R^d)^d$ such that $\grad_\disc u_\disc$ weakly converges in $L^2(\R^d)^d$ to $\bG$ as $h_\disc\to 0$. 
Let us show that $\bG = \grad u$.
Let $\bpsi\in C^\infty_c(\R^d)^d$ be given. Let us consider the term $\terml{ttt}^\disc$ defined by
\[
\termr{ttt}^\disc = \int_{\R^d} \grad_\disc u_\disc(\x) \cdot \bpsi(\x)\d\x.
\]
We get that $\termr{ttt}^\disc =\terml{tt}^\disc + \terml{tu}^\disc$, with
\[
\termr{tt}^\disc = \sum_{K\in\mesh}\sum_{\edge\in\edgescv} \medge (u_\edge - u_K) \ncvedge\cdot \bpsi_K,  
\mbox{ with }\bpsi_K =\frac 1 \mcv\int_K \bpsi(\x)\d\x,
\]
and 
\[
\termr{tu}^\disc = \sum_{K\in\mesh}\sum_{\edge\in\edgescv} R_{K,\edge} u \ \ncvedge   \cdot \int_{D_{K,\edge}}  \bpsi(\x)\d\x.  
\]
We compare $\termr{tt}^\disc$ with $\terml{tt1}^\disc$ defined by
\[
\termr{tt1}^\disc =  \sum_{K\in\mesh}\sum_{\edge\in\edgescv} \medge (u_\edge - u_K) \ncvedge \cdot\bpsi_\edge,
\]
with 
\[
\bpsi_\edge =\frac 1 \medge\int_\edge \bpsi(\x)\d\gamma(\x).
\]
We get that
\[
(\termr{tt}^\disc - \termr{tt1}^\disc)^2 \le \sum_{K\in\mesh}\sum_{\edge\in\edgescv} \frac {\medge} {\dcvedge} (u_\edge - u_K)^2  
\sum_{K\in\mesh}\sum_{\edge\in\edgescv} \medge\dcvedge |\bpsi_K - \bpsi_\edge|^2,
\]
which leads to $\lim\limits_{h_\disc\to 0} (\termr{tt}^\disc - \termr{tt1}^\disc) = 0$. 

Since
\[
\termr{tt1}^\disc =   - \sum_{K\in\mesh}\sum_{\edge\in\edgescv} \medge u_K \ncvedge\cdot \bpsi_\edge = - \int_{\R^d} \Pi_\mesh u_\disc(\x)\div\bpsi(\x)\d\x,
\]
we get that $\lim\limits_{h_\disc\to 0} \termr{tt1}^\disc = - \int_{\R^d} u(\x) \div\bpsi(\x)\d\x$. 
Let us now turn to the study of $\termr{tu}^\disc$. Noting again that \refe{propRnulle} holds, we have: 
\[
\termr{tu}^\disc = \sum_{K\in\mesh}\sum_{\edge\in\edgescv} R_{K,\edge} u \ \ncvedge   \cdot \int_{D_{K,\edge}}  (\bpsi(\x)-\bpsi_K)\d\x.  
\]
Since $\bpsi$ is a regular function, there exists $C_\bpsi$ only depending on $\bpsi$ such that
$\vert \int_{D_{K,\edge}}  (\bpsi(\x)-\bpsi_K)\d\x \vert \le C_\bpsi h_\disc \dfrac{\medge \dcvedge}{d}.$ 
{F}rom \refe{etdeux} and the Cauchy-Schwarz inequality, we thus get:
\[
\lim_{h_\disc\to 0}  \termr{tu}^\disc = 0.
\]

This proves that the function $\bG\in L^2(\R^d)^d$ is a.e. equal to $\grad u$ in $\R^d$. 
Since $u=0$ outside of $\O$, we get that $u\in H^1_0(\O)$, and the uniqueness of the limit implies that the whole family $\grad_\disc u_\disc$ weakly converges in $L^2(\R^d)^d$ to $\grad u$ as $h_\disc\to 0$.

\end{proof}

Note that the proof that $u\in H^1_0(\O)$ also results from \refe{comparnorm}, which allows us to apply Lemma \ref{reglim} of the Appendix in the particular case $p=2$.
Let us also remark that several  discrete gradients could be chosen, which satisfy the weak compactness property (see for instance the proof of Lemma \ref{reglim}).
However, we emphasise that the choice of the specific gradient \eqref{defgradap} also stems from coercivity and consistency issues.
Let us now state the discrete gradient consistency property.

\begin{lemma}[Discrete gradient consistency]\label{consgrad}
Let $\disc$ be a discretisation of $\Omega$ in the sense of Definition \ref{adisc}, and let $\theta \ge \theta_\disc$ be given. 
Then, for any function $\varphi \in C^2(\overline \Omega)$, there exists $\ctel{200}$  only depending on $d$, $\theta$ and $\varphi$ such that:
\begin{equation}
\Vert \nabla_\disc P_\disc \varphi -\nabla \varphi \Vert_{(L^{\infty}(\Omega))^d} \le \cter{200} h_\disc, \end{equation}
where $\nabla_\disc$ is defined by \refe{defgradap}-\refe{defgrad}.
\end{lemma}
\begin{proof}
{F}rom definitions \refe{defgrad} and \refe{defgradKedge}  we get
$$\vert \nabla_{K,\edge} P_\disc \varphi - \nabla \varphi (\x_K) \vert \le \vert  \nabla_{K} P_\disc \varphi - \nabla \varphi (\x_K) \vert + \vert   R_{K,\edge} P_\disc \varphi \vert.$$
{F}rom \refe{defgradap}, we have, for any $K \in \mesh$,  
\begin{eqnarray*}
\nabla_K P_\disc \varphi & =&\dfrac{1}{\mcv } \sum_{\edge \in \edgescv} \medge (\varphi(\xedge) -\varphi(\xcv))\ncvedge \\
 &=&\dfrac{1}{\mcv } \sum_{\edge \in \edgescv} \medge 
\bigl( \nabla \varphi(\xcv) \cdot (\xedge -\xcv) + h_K^2 \rho_{K,\edge} \bigr)\ncvedge,
\end{eqnarray*}
where $\vert \rho_{K,\edge} \vert \le C_\varphi$ with $C_\varphi$ only depending on $\varphi$.
Thanks to  \refe{magical} and to the regularity of the mesh, we get
$$
\vert \nabla_K P_\disc \varphi  -\nabla  \varphi (\xcv) \vert \le \dfrac{1}{\mcv } \sum_{\edge \in \edgescv} \medge  h_K^2 \vert \rho_{K,\edge} \vert \le  h_K \ d \  C_\varphi \theta.
$$
{F}rom this last inequality, using Definition \ref{defrkedge}, we get
\begin{eqnarray*}
\vert   R_{K,\edge} P_\disc \varphi \vert & = & \dfrac{\sqrt{d}}{\dcvedge} \vert \varphi(\xedge) -\varphi(\xcv) - \nabla_K  P_\disc \varphi \cdot (\xedge -\xcv ) \vert\\
 &\le & \dfrac{\sqrt{d }}{\dcvedge} \left( h_K^2 \rho_{K,\edge} +  h_K^2 \ d \  C_\varphi \theta \right)\\
 & \le & \sqrt{d} \theta (h_K  C_\varphi + h_K  d C_\varphi\theta),
\end{eqnarray*}  
 which concludes the proof.
\end{proof}

We now give the abstract properties of the discrete fluxes, which are necessary to prove the convergence of the general scheme \refe{scheme}, \refe{defbilF}, and then prove that the fluxes that we constructed in Section \ref{sec-cons-flux} indeed satisfy these properties.

\begin{definition}[Continuous, coercive, consistent and symmetric families of fluxes]\label{ccsf}~

Let ${\mathcal F}$ be a family of discretisations in the sense of definition \ref{adisc}.
For $\disc = (\mesh, \edges, \centers) \in {\mathcal F}$, $K \in \mesh$ and $\edge \in \edges$, we denote by  
$F_{K,\sigma}^\disc$  a linear  mapping  from $X_{\disc}$ to $\R$, and we denote by $ \Phi= ((F^\disc_{K,\sigma})_{\stackrel{K \in \mesh}{\edge \in \edges}})_{\disc \in \mathcal F}$.  We consider the
 bilinear form defined by
\be
\langle u,v \rangle_F = \sum_{K \in \mesh} \sum_{\edge \in \edgescv}  F^\disc_{K,\sigma}(u)(v_K-v_\edge)\qquad \forall
(u,v) \in  X_{\disc}^2.
 \label{deffobil}\ee

The family of numerical fluxes $ \Phi$ is said to be continuous if there exists $M>0$  such that
\begin{equation}
\langle u,v \rangle_F \le M \vert u \vert_{X} \vert v \vert_{X}\qquad \forall (u,v) \in X_{\disc}^2,\ 
\forall \disc=(\mesh,\edges,\centers)\in  \mathcal F. 
\label{fluxcontinuous}
\end{equation}

The family of numerical fluxes $ \Phi$ is said to be coercive if there exists $\alpha >0$  such that
\begin{equation}
 \alpha \vert u \vert^2_{X} \le 
\langle u,u \rangle_F\qquad \forall u \in X_{\disc} \  \forall \disc=(\mesh,\edges,\centers)\in  \mathcal F. 
\label{fluxcoercif}
\end{equation}

%A family of such linear mappings 
The family of numerical fluxes $ \Phi$ is said to be consistent (with Problem \refe{eq1}--\refe{cl0}) if for any family $(u_\disc)_{\disc\in\mathcal F}$ satisfying: 
\begin{itemize}
\item $u_\disc\in X_{\disc,0}$ for all $\disc\in\mathcal F$,
\item there exists $C>0$ with $\vert u_\disc\vert_X \le C$ for all $\disc\in\mathcal F$,
\item there exists $u\in L^2(\O)$ with $\lim\limits_{h_\disc \to 0}\Vert \Pi_\mesh u_\disc - u\Vert_{L^2(\O)} = 0$
(recall that, from Lemma \ref{reglim}, we get that $u\in H^1_0(\O)$), 
\end{itemize}
 then
\be
\lim_{h_\disc \to 0} \langle u_\disc,P_{\disc} \varphi\rangle_F = 
\int_\O \Lambda(\x) \nabla \varphi(\x) \cdot \nabla u(\x) \d\x \qquad\forall \varphi \in C^\infty_c(\O).
\label{consisflux}
\ee

Finally the family of numerical fluxes  $\Phi$  is said to be  symmetric if 
\[
\langle u,v \rangle_F = \langle v,u \rangle_F\qquad \forall (u,v) \in  X_{\disc}^2,\  \forall \disc=(\mesh,\edges,\centers)\in  \mathcal F.
\]

\end{definition}

We now show that the family of fluxes defined by  \refe{defy}-\refe{defflups}  satisfies the definition of a consistent, coercive and symmetric family of fluxes. 
Recall that the SUSHI scheme \refe{schemepratique} is studied numerically in Section \ref{secresnum}  with this choice for the family of fluxes.

\begin{lemma}[Flux properties]\label{casparthyb}
Let ${\mathcal F}$ be a family of discretisations in the sense of Definition \ref{adisc}.
We assume that there exists $\theta >0$ with 
\be\ba
\theta_{\disc}\le \theta\qquad \forall \disc=(\mesh,\edges,\centers)\in {\mathcal F},
\ea\label{regultheta}\ee
where $\theta_{\disc}$ is defined by \refe{defetadisc}. 
Let  $\Phi=((F^\disc_{K,\sigma})_{\stackrel{K \in \mesh}{\edge \in \edgescv}})_{\disc \in \mathcal F}$ be the family of fluxes defined by \refe{defy}-\refe{defflups}. 
Then, the family $\Phi$ is a  continuous, coercive, consistent and symmetric family of numerical fluxes in the sense of Definition \ref{ccsf}.
\end{lemma}

\begin{proof} Since the family of fluxes is defined by \refe{defy}-\refe{defflups}, it satisfies \refe{firstidea}, and therefore we have:
\[
\langle u,v\rangle_F = \int_\O \grad_\disc u(\x)\cdot \Lambda(\x)\grad_\disc v(\x)  \d\x \qquad \forall u,v\in X_{\disc}.
\]
Hence the property $\langle u,v\rangle_F = \langle v,u\rangle_F$ holds.
The continuity and coercivity of the family $\Phi$ result   from Lemma \ref{lemgrad} and the properties of $\Lambda$, which give: $\langle u,v \rangle_F \le \overline{\lambda} \Vert \nabla_\disc u \Vert_{L^2(\O)}\Vert \nabla_\disc v \Vert_{L^2(\O)}$ and 
$\langle u,u\rangle_F \ge \underline{\lambda} \Vert \nabla_\disc u \Vert_{L^2(\O)}^2$ for any $u,v \in X_\disc.$   
The consistency results from the weak and strong convergence properties in Lemmas \ref{consflux} and \ref{consgrad}, which give $\nabla_\disc u_\disc  \to \nabla u$ weakly in $L^2(\O)$ and $\nabla_\disc P_\disc \varphi \to \nabla \varphi$ in $L^2(\O)$ as the mesh size tends to 0.
\end{proof}

\begin{theorem}[Convergence]\label{cvgce}
Let ${\mathcal F}$ be a family of discretisations in the sense of Definition \ref{adisc}, for any $\disc \in {\mathcal F}$,  let  $\bary \subset \edgesint$ and $(\beta_\edge^K)_{\stackrel{K\in\mesh}{\edge\in\edgesint}}$  satisfying \refe{propbary}. 
Assume that there exists  $\theta >0 $  such that $\theta_{\disc,\bary} \le \theta$,  for all $\disc \in {\mathcal F}$, where $\theta_{\disc,\bary}$ is defined by \refe{regutousazimut}. 
Let $\Phi = ((F^\disc_{K,\sigma})_{\stackrel{K \in \mesh}{\edge \in \edges}})_{\disc \in \mathcal F}$ be a continuous,  coercive and symmetric and consistent family of numerical fluxes in the sense of Definition \ref{ccsf}. 
Let $(u_\disc)_{\disc\in\mathcal F}$  be the family of functions  satisfying \refe{scheme} for all
$\disc\in\mathcal F$.
Then $\Pi_\mesh u_\disc$ converges in $L^2(\O)$ to the unique solution $u$ of \refe{ellgenf} as $h_\disc\to 0$.
Moreover $\grad_\disc u_\disc$ converges to $\grad u$ in $L^2(\O)^d$ as $h_\disc\to 0$.
\end{theorem}

\begin{proof}
Letting $v = u_\disc$ in \refe{scheme} and applying the Cauchy-Schwarz inequality yields
\[
\langle u_\disc, u_\disc\rangle_F = \int_\O f(\x)  \Pi_\mesh u_\disc(\x) \d\x \le \Vert f\Vert_{L^2(\O)}\Vert \Pi_\mesh u_\disc\Vert_{L^2(\O)}.
\]
We apply the Sobolev inequality \refe{insob} with $p=2$, which gives in this case  
$$\Vert \Pi_\mesh u_\disc\Vert_{L^2(\O)} \le \ctel{csobjan} \Vert \Pi_\mesh u \Vert_{1,2,\mesh}.$$
Using \refe{comparnorm} and the consistency of the family $\Phi$ of fluxes, we then have
\[
\alpha \vert \Pi_\mesh u_\disc\vert_X^2 \le \cter{csobjan}\Vert f \Vert_{L^2(\O)} \vert u_\disc\vert_X.
\]
This  leads to the inequality
\be
\Vert u_\disc \Vert_{1,2,\mesh}\le \vert u_\disc\vert_X\le \frac{\cter{csobjan}} {\alpha} \Vert f\Vert_{L^2(\O)}.
\label{estimu}\ee
Thanks to Lemma \ref{reglim}, we get the existence of $u\in H^1_0(\O)$, and of a subfamily extracted from
$\mathcal F$, such that $\Vert\Pi_\mesh  u_\disc - u\Vert_{L^2(\O)}$ tends to $0$ as $h_\disc\to 0$. 
For a given $\varphi\in C^\infty_c(\O)$, let us take $v = P_{\disc,\bary} \varphi$ in  \refe{scheme} (recall that $P_{\disc,\bary} \varphi \in X_{\disc,\bary}$).
We get
\[
\langle u_\disc,P_{\disc,\bary} \varphi\rangle_F = \int_\O f(\x) P_\mesh\varphi(\x) \d\x.
\]
Let us remark that, thanks to the continuity of the family $\Phi$ of fluxes, we have
\[
\langle u_\disc,P_{\disc,\bary} \varphi - P_{\disc} \varphi\rangle_F \le  
M \frac{\cter{afdinsobp}} {\alpha} \Vert f\Vert_{L^2(\O)}\ \vert P_{\disc,\bary} \varphi - P_{\disc} \varphi\vert_X.
\]
Thanks to \refe{propbary} and \refe{regutousazimut}, we get the existence of $C_\varphi$ only depending on $\varphi$ (through its second order partial derivatives) such that, for all $K\in\mesh$ and all $\edge\in\bary\cap\edgescv$,
\be
\vert \sum_{L\in\mesh}\beta_\edge^L \varphi(\xcvv) - \varphi(\xedge)\vert \le  
 \sum_{L\in\mesh} \vert\beta_\edge^L\vert |\xcvv - \xedge|^2 C_\varphi \le 
\theta_{\disc,\bary} C_\varphi h_K^2.
\label{ineqphibar}\ee
We can then  deduce
\be
\lim_{h_\disc\to 0} \vert P_{\disc,\bary} \varphi - P_{\disc} \varphi\vert_X = 0.
\label{diffpdbp}\ee
Thanks to the $\mathcal F$-extracted subfamily properties, we may apply the consistency hypothesis on the family $\Phi$ of fluxes, which gives
\[
\lim_{h_\disc\to 0} \langle u_\disc,P_{\disc} \varphi\rangle_F  = 
\int_\O \Lambda(\x) \nabla \varphi(\x) \cdot \nabla u(\x) \d\x.
\]
Gathering the two results above leads to
\[
\lim_{h_\disc\to 0} \langle u_\disc,P_{\disc,\bary} \varphi\rangle_F  = 
\int_\O \Lambda(\x) \nabla \varphi(\x) \cdot \nabla u(\x) \d\x,
\]
which concludes the proof of the following equality
\[
\int_\O \Lambda(\x) \nabla \varphi(\x) \cdot \nabla u(\x) \d\x = \int_\O f(\x)  \varphi(\x) \d\x.
\]
Therefore, $u$ is the unique solution of \refe{ellgenf}, and 
we get that the whole family $(u_\disc)_{\disc\in\mathcal F}$ converges
to $u$ as $h_\disc\to 0$.

\smallskip

Let us now prove the second part of the theorem.

Let $\varphi\in {\rm C}^\infty_c(\O)$ be given (this function is devoted to approximate $u$ in $H^1_0(\O)$). 
Thanks to the Cauchy-Schwarz inequality, we have
\[
\int_\O |\grad_{\disc}u_\disc(\x) - \grad u(\x)|^2 \d \x \le 3 \ (\terml{2}^\disc+\terml{6}^\disc+\terml{5}),
\]
with $\termr{2}^\disc = \int_\O \vert \grad_{\disc}u_\disc(\x) - \grad_{\disc}P_{\disc}\varphi(\x) \vert^2 \d \x, $ $\termr{6}^\disc =\int_\O  \vert \grad_{\disc}P_{\disc}\varphi(\x) - \grad\varphi(\x) \vert^2 \d \x,$ 
 and $
\termr{5} =\int_\O  \vert \grad\varphi(\x) - \grad  u (\x) \vert^2 \d \x.
$
Thanks to Lemma \ref{consgrad},  we have
$
\lim_{h_\disc \to 0} \termr{6}^\disc = 0.
$

Thanks  to Lemma \ref{lemgrad} and to the coercivity of the family of fluxes, there exists $\ctel{bG}$ such that
\[
\Vert \grad_{\disc} v \Vert_{L^2(\O)^d}^2 \le \cter{102}^2\vert v\vert_X^2 \le   \cter{bG}\langle v,v\rangle_F \qquad \forall v\in X_{\disc},
\]
with $\cter{bG} = \frac{\cter{102}^2} {\alpha}$.
Taking $v = u_\disc - P_{\disc}\varphi$, we have
\[
\termr{2}^\disc \le \cter{bG} (\langle u_\disc,u_\disc\rangle_{F}
- 2 \langle u_\disc,P_{\disc}\varphi\rangle_F + 
\langle P_{\disc}\varphi,P_{\disc}\varphi\rangle_F ).
\]
By Theorem \ref{cvgce} and thanks to  and consistency of the family of fluxes, we get
$$
\lim_{h_\disc \to 0} \langle u_\disc,P_{\disc}\varphi\rangle_F =  \int_\O \grad  u (\x)\cdot\Lambda(\x) \grad\varphi(\x) \d \x. 
%\label{vvv31}\ee
$$
The sequence $\vert P_{\disc}\varphi \vert_X$ is bounded; using the regularity of $\varphi$, the regularity hypotheses of the family of discretisations, together with the consistency of the family of fluxes implies that
$$
\lim_{h_\disc \to 0} \langle P_{\disc}\varphi,P_{\disc}\varphi\rangle_F  = \int_\O \grad\varphi(\x)\cdot\Lambda(\x) \grad\varphi(\x) \d \x.
$$
Remarking that passing to the limit $h_\disc \to 0$ in \refe{scheme} with $v=u_\disc$ provides that $\langle u_\disc,u_\disc \rangle_F$ converges to $\int_\Omega \nabla u\cdot\Lambda \nabla u \d \x$, we get that
\[
\lim_{h_\disc \to 0} \langle u_\disc - P_{\disc}\varphi,u_\disc - P_{\disc}\varphi\rangle_{F} =
\int_\Omega \nabla (u-\varphi)\cdot\Lambda \nabla  (u-\varphi) \d \x \le \overline \lambda
\int_\O \vert \grad  u - \grad\varphi\vert^2 \d \x,
\]
which yields
\[
\limsup_{h_\disc \to 0}\termr{2}^\disc \le \cter{bG} \overline \lambda
\int_\O \vert\grad  u - \grad\varphi\vert^2 \d \x.
\]
{F}rom the above results, we obtain that there exists $\ctel{toto}$, independent of $\disc$, such that
\[
\int_\O \vert\grad_{\disc}u_\disc(\x) - \grad  u (\x)\vert^2 \d \x\le \cter{toto} \int_\O \vert\grad\varphi(\x) - \grad u (\x)\vert^2 \d \x + \terml{7}^\disc,
\]
with (noting that $\varphi$ is fixed)  $\lim_{h_\disc \to 0} \termr{7}^\disc = 0.$
Let $\varepsilon >0$;  we may choose $\varphi$ such that $\int_\O \vert \grad\varphi(\x) - \grad u (\x)\vert^2 \d \x \le \varepsilon$, and we may then choose $h_\disc $ small enough so that $\termr{7}^\disc\le \varepsilon$. 
This completes the proof that
\be
\lim_{h_\disc \to 0} \int_\O \vert\grad_{\disc}u_\disc(\x) - \grad  u (\x)\vert^2 \d \x = 0
\label{vvv5}\ee
in the case of a general continuous, coercive, consistent and symmetric family of fluxes.
\end{proof}

Let us write an error estimate in the particular case $\Lambda = \rm{Id}$, assuming a regular exact solution to \eqref{ellgenf}.
\begin{theorem}[Error estimate, isotropic case]\label{eresthyp}
We consider the particular case $\Lambda = \rm{Id}$, and we assume that the solution $u\in H^1_0(\O)$ of \refe{ellgenf} is in  $C^2(\overline{\O})$. Let $\disc=(\mesh,\edges,\centers)$ be a discretisation in the sense of Definition \ref{adisc}, let $\bary\subset\edgesint$ be given, let $\bary =   (\beta_\edge^K)_{\edge\in\bary,K\in\mesh}\subset \R$ such that \refe{propbary} holds, and let $\theta \ge \theta_{\disc,\bary}$ be given 
(see \refe{regutousazimut}).
% re ici pas de family de continuous etc fluxes car une seule discretisation !!!
Let $(F_{K,\sigma})_{K \in \mesh,\edge \in \edges}$  be a family of linear  mappings  from $X_{\disc}$ to $\R$, such that there exists $\alpha >0$ with
\begin{equation}
 \alpha \vert v \vert^2_{X} \le  \langle v,v \rangle_F\qquad \forall v \in X_{\disc}, 
\label{fluxcoercifester}\end{equation}
defining $\langle \cdot,\cdot \rangle_F$ by \refe{deffobil}.
%%%%%%%%%%% re 
We denote by
\begin{equation}
E(u) = \left(\sum_{K\in \mesh}\sum_{\edge \in \edgescv}\frac {\dcvedge} {\medge} \left(  F_{K,\sigma}(P_{\disc,\bary}  u) +
\int_\edge  \nabla u(\x) \cdot \ncvedge \d \gamma(\x) \right)^2\right)^{1/2}.
\label{consisfluxdis}
\end{equation}
Then the solution $u_\disc$ of \refe{scheme} satisfies that there exists $\ctel{esterhyb}$, only depends on $\alpha$ and on $\theta$,
such that
\be
\Vert \Pi_\mesh u_\disc - P_\mesh u\Vert_{L^2(\O)} \le \cter{esterhyb} E(u),
\label{erresthybeq}
\ee
and satisfies that there exists $\ctel{esterhybg}$, only depending on  $\alpha$, $\theta$ and $u$ such that
\be
\Vert \grad_\disc u_\disc - \grad u\Vert_{L^2(\O)^d} \le \cter{esterhybg} \left(E(u) + h_\disc\right).
\label{erresthybeqgh}\ee
Moreover, in the particular case where $(F_{K,\sigma})_{K \in \mesh,\edge \in \edges}$ is defined by \refe{defy}-\refe{defflups},
there exists  $\ctel{caspartic}$, only depending on $\alpha$, $\theta$ and $u$, such that
\begin{equation}
E(u)\le \cter{caspartic} h_\disc.
\label{consisfluxdish}
\end{equation} 
\end{theorem}
\begin{remark}[Extensions of the error estimate]
Note also that the extension of Theorem \ref{eresthyp} to the case $u\in H^2(\O)$ is possible for $d=2$ or $d=3$.
However it would demand  a rather longer and more technical proof and is not expected to provide more information on the link between accuracy and the regularity of the mesh than the result presented here.
In the case of the pure hybrid scheme (HFV, $\bary = \emptyset$), an error estimate could however be obtained by assuming $u$ piece-wise to be $H^2$. 
Such error estimates were also obtained for pure hybrid schemes of the mimetic type by using the tools of the mixed finite element theory (see e.g. \cite{bre-05-fam}). 
If $\bary \not = \emptyset$, one must furthermore assume that the barycentric formulae \eqref{ecrbar}-\eqref{propbary} or  \eqref{ecrbarbis}-\eqref{propbarybis} are written with unknowns located in the same regularity zone, as explained in Remark \ref{rembarycons}.
Nevertheless such error estimates are not possible for general $L^\infty$ diffusion operators, since in such a case the maximal regularity of the continuous solution is $H^1_0(\O)$. 
Then, by interpolation, one may get some error estimates if the continuous solution is in $H^{1}_0(\O)\cap H^{s}(\O)$ as in the classical finite element framework.
%Our concern is not to focus on such results, but to identify the origin of the error terms in a regular case. %(the following result could easily be extended to cases with piecewise constant diffusion operators). 
\end{remark}
\begin{proof}
%{F}rom \refe{ellgenf}, letting $v = P_{\disc,\bary} u - u_\disc \in X_{\disc,\bary}$, we get that
Let $v \in X_\disc$, since $-\Delta u = f$, we get:
\begin{equation}
-\sum_{K\in\mesh} v_K \int_K \lap u(\x)\d\x = \int_\O f(\x) \Pi_\mesh v(\x)\d\x.
\label{ilaresondraler}
\end{equation}
Thanks to the following equality (recall that $u \in C^2(\overline{\Omega})$ and therefore $\nabla u \cdot \ncvedge $ is defined on each edge $\edge$)
\[
-\sum_{K\in\mesh} v_K \int_K \lap u(\x)\d\x = -\sum_{K\in \mesh}\sum_{\edge \in \edgescv}  (v_K - v_\edge)\int_\edge  \nabla u(\x) \cdot \ncvedge \d \gamma(\x),
\]
we get that
\[
\langle P_{\disc,\bary} u, v\rangle_F = \int_\O f(\x)\Pi_\mesh v(\x)\d \x + 
\sum_{K\in \mesh}\sum_{\edge \in \edgescv} \left(F^\disc_{K,\sigma}(P_{\disc,\bary} u) + 
\int_\edge  \nabla u(\x) \cdot \ncvedge \d \gamma(\x) \right) (v_K - v_\edge).
\]
Taking $v = P_{\disc,\bary} u- u_\disc \in X_{\disc,\bary}$ in this latter equality and using \refe{ilaresondraler} we get 
\[
\langle v  , v\rangle_F = \sum_{K\in \mesh}\sum_{\edge \in \edgescv} \left(F^\disc_{K,\sigma}(P_{\disc,\bary} u) + 
\int_\edge  \nabla u(\x) \cdot \ncvedge \d \gamma(\x) \right) (v_K - v_\edge),
\]
which leads, using \refe{fluxcoercifester} and the Cauchy-Schwarz inequality, to
\be
\alpha \vert v \vert_{X} \le E(u).
\label{errvh}\ee
Using \refe{comparnorm} and  the Sobolev inequality \refe{insob} with $p=2$ provides the conclusion of \refe{erresthybeq}.
Let us now prove \refe{erresthybeqgh}. We have
\[
\Vert \grad_\disc u_\disc - \grad u\Vert_{L^2(\O)^d}\le 
\Vert \grad_\disc u_\disc - \grad_\disc P_{\disc,\bary} u\Vert_{L^2(\O)^d}+
\Vert \grad_\disc P_{\disc,\bary} u- \grad u\Vert_{L^2(\O)^d}.
\]
The bound of the first term in the above right-hand side is bounded thanks to Lemma \ref{lemgrad} and \refe{errvh}. 
The inequality 
$\Vert \grad_\disc P_{\disc,\bary} u- \grad u\Vert_{L^2(\O)^d} \le \ctel{cnouv} h_\disc $ 
is obtained thanks to Lemma \ref{consgrad} and using a similar inequality to \refe{ineqphibar}, replacing $\varphi$ by $u$.

\medskip

Let us now turn to the proof of \refe{consisfluxdish} in the particular case
 where the family of fluxes is defined by \refe{defy}-\refe{defflups}. Indeed, we get in this case
that, for all $v\in X_\disc$,
\[
F_{K,\sigma}(v) = -\sum_{\edgep\in\edgescv} (\grad_K v + R_{K,\edgep} v\ \ncvedgep)\cdot \frac {\medgep\dcvedgep} {d} 
\mathbi{y}^{\edgep \edge},
\]
with
\[
\mathbi{y}^{\edgep  \edge}= \left\{\begin{array}{ll}
\dsp \frac \medge \mcv \ncvedge +  \frac {\sqrt{d}} {\dcvedge} \left(1 - \frac\medge  \mcv \ncvedge \cdot (\xedge - \xcv)\right)\ncvedge & \mbox{ if } \edge = \edgep \\
\dsp \frac \medge \mcv \ncvedge - \frac {\sqrt{d}} {\dcvedgep \mcv}\medge \ncvedge \cdot (\xedgep - \xcv)\ncvedgep & \mbox{ otherwise }.
\ea\right.\]
Using \refe{magical}, we get that 
\[
\sum_{\edgep\in\edgescv} \frac {\medgep\dcvedgep} {d}  \mathbi{y}^{\edgep \edge} = \medge \ncvedge.
\]
Since there exists  $\ctel{truc}\in \R_+$ such that  $ |R_{K,\edgep} P_{\disc,\bary} u| \le \cter{truc} h_K$,  there exists some $\ctel{truc2}\in \R_+$ with
\[
\left\vert F_{K,\sigma}(P_{\disc,\bary}  u) +
\int_\edge  \nabla u(\x) \cdot \ncvedge \d \gamma(\x)\right\vert \le \cter{truc2} \medge h_K.
\]
This leads to the conclusion of \refe{consisfluxdish}.
\end{proof}

\section{Discrete functional analysis}\label{secsobdis}

This section is devoted to some results of functional analysis that are useful for the proof of convergence of numerical schemes when the approximate solution is  piece-wise constant on the mesh.
Although some of the results presented here were already introduced in previous works of the authors, they were mostly presented (even when not needed, see \cite[Remark 9.13 p. 793]{book}) in the framework of ``admissible" meshes, that is meshes with an orthogonality condition. 

We recall that in the proof of the main convergence Theorem \ref{cvgce}, we first obtain from the scheme some estimates on the approximate solutions in the discrete $H^1$ norm. 
We now show how, from a general discrete $W^{1,p}$ estimate (this generalisation to $p \ne 2$ is useful in the case of nonlinear problems) we obtain a discrete $L^q$ estimate for some $q >p$ (Lemma \ref{sob}). 
We then obtain a certain compactness result in $L^1$ (Lemma \ref{kolmun} and therefore in $L^p$ (Lemma \ref{kolmp}), which in turn allows to show that the limit of the approximate solution is in $W^{1,p}_0(\Omega)$ (Lemma \ref{reglim}). 

\subsection{Discrete Sobolev embeddings}

\subsubsection{Discrete embedding of $W^{1,1}$ in $L^{1^\star}$}

The discrete Sobolev embedding  of $W^{1,1}$ in $L^{1^\star}$ requires less assumptions on the mesh than those given  in Definition~\ref{adisc}.
We therefore introduce a larger class of meshes in the following definition.

\begin{definition}[Polyhedral partition of $\O$]\label{defpart}
Let $d \ge 1$ and let $\O$ be an open bounded set in $\R^d$, whose boundary is a finite union of part of hyperplanes. 
A polyhedral partition  $\mesh$  of $\O$ is a finite partition of $\O$ such that each element $K$ of this partition is measurable and has a boundary  $\partial K$ that is composed of a finite union of parts of hyperplanes (the facets of $K$) denoted by $\edge$: $\partial K = \cup_{\sigma \in \edgescv}  \sigma$.
%Let $\mesh$ be such a partition and
Let $\edges$ be the set of the facets of all the elements of $\mesh$:  $\edges =  \cup_{K\in\mesh} \edgescv$.
If $\sigma \in \edges$ is a facet of this partition, one denotes by  $\medge$ the $(d-1)$--Lebesgue measure of $\sigma$. 
Let $H_\mesh(\O)$ be the set of functions from $\O$ to $\R$, constant on each element of $\mesh$.
Let $u \in H_\mesh(\O)$. 
If $\sigma \in \edgescv \cap \edgescvv$ (that is $\sigma$ is a facet such that $\sigma \subset\overline K \cap \overline L$),  one sets $D_\sigma u=\vert u_K-u_L \vert$.
If $\sigma \in \edges$  is on the boundary of $\O$ and $K \in \mesh$ (that is $\sigma= \partial \O \cap \overline K$), one sets $D_\sigma u=\vert u_K \vert$. For $u \in H_\mesh(\O)$, one sets
\be
\norm{u}{1,1,\mesh}=\sum_{\sigma \in \edges } \medge D_\sigma u.
\label{ppo}
\ee
\end{definition}

\begin{lemma} Let $d \ge 1$ and let $\O$ be an open bounded set of $\R^d$, whose boundary is a finite union of parts of hyperplanes. 
Let $\mesh$ be a polyhedral partition of $\O$ in the sense of Definition \ref{defpart}. 
Then, with the notations of Definition \ref{defpart},
\begin{equation}
\norm{u}{L^{1^\star}(\O)} \le \frac 1 {2\sqrt{d}} \norm{u}{1,1,\mesh} \qquad\forall u \in H_\mesh(\O), \mbox{ with }1^\star =\frac d {d-1}.
\label{insobun}
\end{equation}
\label{sobun}
\end{lemma}
\begin{proof}
Different proofs of this lemma are possible. 
A first proof consists in adapting to this discrete setting the classical proof of the Sobolev embedding due to L. Nirenberg (actually, it gives $1/2$ instead of $1/(2\sqrt{d})$ in \refe{insobun}): it is based  on an induction on $d$. 
This proof is essentially given in \cite[Lemma 9.5 page 790]{book}, with slightly less general hypotheses; in fact the so called orthogonality assumption is not used in  the proof of Lemma 9.5 of \cite{book}.
An easy adaptation of this proof leads to the present lemma (with $1/2$ instead of $1/(2\sqrt{d})$ in \refe{insobun}).

The present proof  makes direct use of  L. Nirenberg's result, namely:
\be
\norm{u}{L^{1^\star}(\R^d)} \le \frac 1 {2d} \norm{u}{W^{1,1}(\R^d)}\qquad\forall u \in W^{1,1}(\R^d),
\label{nirenb}
\ee
where $\norm{u}{W^{1,1}(\R^d)}= \sum_{i=1}^d \norm{D_iu}{L^1(\R^d)}$ and $D_iu$ is the weak
derivative (or derivative in the sense of distributions) of $u$ in the direction $x_i$ (with
$\x=(x_1,\ldots,x_d) \in \R^d$).

\vspace{0.2cm}

For $u \in L^1(\R^d)$, one sets $\norm{u}{BV}=\sum_{i=1}^d \norm{D_iu}{M}$ with, for $i=1,\ldots,d$, $ \norm{D_iu}{M}=\sup\{ \int u \frac {\partial \varphi}{\partial x_i}  \d\x$, $\varphi  \in C^\infty_c(\R^d)$, $\norm{\varphi}{L^\infty(\R^d)} \le 1\}$. 
The function $u$ belongs to the space $BV$ if $u \in L^1(\R^d)$ and $\norm{u}{BV} < \infty$. 
We first remark that \refe{nirenb} is true with $\norm{u}{BV}$ instead of $\norm{u}{W^{1,1}(\R^d)},$ and if $u \in BV$ instead of $W^{1,1}(\R^d)$.
Indeed, to prove this result (which is classical), let $\rho \in C^\infty_c(\R^d,\R_+)$ with $\int \rho \d\x=1$. For $n \in \N^\star$, define $\rho_n=n^d\rho(n \cdot)$.
Let $u \in BV$ and $u_n = u \star \rho_n$ so that, with \refe{nirenb}:
\be
\norm{u_n}{L^{1^\star}(\R^d)} \le \frac 1 {2d} \sum_{i=1}^d \norm{D_i u_n}{L^1(\R^d)}.
\label{nirenbt}
\ee
Since $u_n$ is regular, $\norm{D_i u_n}{L^1(\R^d)}=\norm{D_i u_n}{M}$, and, for $\varphi \in C^\infty_c(\R^d)$, using Fubini's theorem:
$$
\int_{\R^d} u_n \frac {\partial \varphi} {\partial x_i}\d\x = \int_{\R^d} u \frac {\partial } {\partial x_i} (\varphi \star \rho_n) \d\x \le \norm{D_iu}{M} \norm{\varphi}{L^\infty(\R^d)}.
$$
This leads to $\norm{D_i u_n}{L^1(\R^d)} \le \norm{D_iu}{M}$. Since $u_n \tends u$ a.e., as $n \tends \infty$, at least for a sub-sequence, Fatou's lemma gives, from \refe{nirenbt}:
\be
\norm{u}{L^{1^\star}(\R^d)} \le \frac 1 {2d} \norm{u}{BV} \qquad \forall u \in BV.
\label{nirenbtt}
\ee
Let $u \in H_\mesh(\O)$. One sets $u=0$ outside $\O$ so that $u \in L^1(\R^d)$.
One has
$\norm{u}{BV} = \sup\{ \int_{\R^d} u\ \div \varphi \ \d\x$, $\varphi 
\in C^\infty_c(\R^d, \R^d)$, $\norm{\varphi}{L^\infty(\R^d)} \le 1\}$, with
$\norm{\varphi}{L^\infty(\R^d)} = \sup_{i=1,\ldots, d} \norm{\varphi_i}{L^\infty(\R^d)}$ and
$\varphi=(\varphi_1,\ldots,\varphi_d)$.
But, for $\varphi \in C^\infty_c(\R^d, \R^d)$ such that $\norm{\varphi}{L^\infty(\R^d)} \le 1$, an integration by parts on each element of $\mesh$ gives (where ${\bf n}_\sigma$ is a normal vector to $\sigma$ and $\gamma$ is the $(d-1)-$Lebesgue measure on $\sigma$):
$$
 \int_{\R^d} u \ \div \varphi \ \d\x = \sum_{\sigma \in \edges} D_\sigma u \int_\sigma
\vert \varphi \cdot {\bf n}_\sigma \vert d\gamma(\x)\le \sqrt{d} \norm{u}{1,1,\mesh}.
$$
Then, one has $ \norm{u}{BV} \le \sqrt{d}\norm{u}{1,1,\mesh}$ and \refe{nirenbtt} leads
to \refe{insobun}.

%remark for lecteur atentif... le cote embetant de la 1ere proof est cache
%dans l'integration par parties...
\end{proof}

\subsubsection{Discrete embedding of $W^{1,p}$ in $L^{p^\star}$, $1<p<d$}

We now prove a discrete Sobolev embedding for $1<p<d$ and for meshes in the sense of Definition~\ref{adisc}.

\begin{lemma}Let $d > 1$, $1<p<d$ and let $\O$ be a polyhedral open bounded connected subset of $\R^d$.
Let $\disc$ be a discretization on $\O$ in the sense of Definition~\ref{adisc}.
Let $\eta >0$ be such that $\eta \le d_{K,\sigma}/d_{L,\sigma} \le 1/\eta$ for all $\sigma \in \edges$, where $\mesh_\sigma=\{K,L\}$.
%ndrh vire arce qu'on ne voit pas a coui ca sert
%(see the definitions in Section  \ref{cvstudy} and Definition \ref{defpart}). 
Then, there exists $\ctel{afdinsobp}$, only depending on $d$, $p$ and $\eta$ such that
\be
\norm{u}{L^{p^\star}(\O)} \le \cter{afdinsobp} \norm{u}{1,p,\mesh} \qquad  \forall u \in H_\disc(\O),
\label{insobp}
\ee
where $p^\star= \frac {pd}{d-p}$ and
\be
\norm{u}{1,p,\mesh}^p=\sum_{K \in \mesh} \sum_{\sigma \in \edgescv} \medge d_{K,\sigma}
\left(\frac {D_\sigma u} {d_{\sigma}}\right)^p,
\label{defnorunpd}\ee
with $d_{\sigma}=d_{K,\sigma} + d_{L,\sigma}$, if $\mesh_{\sigma}=\{K,L\}$, and
$d_\sigma=d_{K,\sigma}$, if $\mesh_\sigma=\{K\}$.
\label{sobp}
\end{lemma}
\begin{proof}
We follow here  L. Nirenberg's proof of the Sobolev embedding. 
Let $\alpha$ be such that $\alpha 1^\star= p^\star$ (that is $\alpha=p(d-1)/(d-p) >1$).
Let $u \in H_\disc(\O)$. 
Inequality \refe{insobun} applied with $\vert u \vert^\alpha$ instead of $u$ leads to:
$$
\left(\int_\O \vert u \vert^{p^\star} \d\x\right)^{\frac {d-1}{d}} \le
\sum_{\sigma \in \edges} \medge D_\sigma \vert u \vert^\alpha.
$$
For $\sigma \in \edgesint$, $\mesh_\sigma = \{K,L\}$, one has $D_\sigma \vert u \vert^\alpha
\le \alpha (\vert u_K \vert^{\alpha -1} + \vert u_L \vert^{\alpha -1}) D_\sigma u$. 
For $\sigma \in \edgesext$, $\mesh_\sigma = \{K\}$, one has $D_\sigma \vert u \vert^\alpha
\le \alpha \vert u_K \vert^{\alpha -1} D_\sigma u$. 
This yields:
\be
\left(\int_\O \vert u \vert^{p^\star} \d\x\right)^{\frac {d-1}{d}} 
\le \sum_{K \in \mesh} \sum_{\sigma \in \edgescv} \medge \alpha\vert u_K \vert^{\alpha -1}
D_\sigma u,
\label{newpp}
\ee
For all $\sigma \in \edges$, one has $1 \le \frac {1+\eta} {\eta}\frac {d_{K,\sigma}} {d_\sigma}$,
if  $\sigma \in \edgesint$, $\mesh_\sigma = \{K,L\}$, or if $\sigma \in \edgesext$, $\mesh_\sigma = \{K\}$.
Then, H\"older's inequality applied to \refe{newpp} yields, with $q=p/(p-1)$:
\be
(\int_\O \vert u \vert^{p^\star} \d\x)^{\frac {d-1}{d}} \le \alpha \frac {1+\eta} {\eta}
(\sum_{K \in \mesh} \sum_{\sigma \in \edgescv}
\medge d_{K,\sigma} \vert u_K \vert^{(\alpha-1)q})^{\frac 1 q} \norm{u}{1,p,\mesh}.
\label{insobpp}
\ee
Since $(\alpha-1)q=p^\star$, one has:
$$
\sum_{K \in \mesh} \sum_{\sigma \in \edgescv}
\medge d_{K,\sigma} \vert u_K \vert^{(\alpha-1)q}=d \int_\O \vert u \vert^{p^\star} \d\x.
$$
Then, noticing that $(d-1)/d - 1/q = 1/p^{\star}$, we deduce \refe{insobp} follows from \refe{insobpp} with $\cter{afdinsobp}=\alpha \frac {1+\eta} {\eta}d^{1/q}$ only depending on $d$, $p$ and $\eta$.
\end{proof}

\subsubsection{Discrete embedding of $W^{1,p}$ in $L^{q}$, for some $q>p$}
Let $1 \le p < \infty$, we now deduce from Lemma \ref{sob} the following lemma, which gives the discrete  embedding of $W^{1,p}$ in $L^{q}$, for some $q>p$.
 \begin{lemma}Let $d \ge  1$, $1\le p< \infty$ and  let $\O$ be a polyhedral open bounded connected subset of $\R^d$.
Let $\disc$ be a mesh of $\O$ in the sense of Definition~\ref{adisc}.
Let $\eta >0$ be such that $\eta \le d_{K,\sigma}/d_{L,\sigma} \le 1/\eta$ for all $\sigma \in \edges$, where $\mesh_\sigma=\{K,L\}$. 
Then, there exists $q>p$ only depending on $p$ and there exists $\ctel{afdinsob}$, only depending on $d$, $\O$, $p$ and $\eta$ such that 
%ndrh vire parce qu'on ne sait pas a auoi ca sert (the definitions in Section \ref{cvstudy} and Definition \ref{defpart})
\be
\norm{u}{L^{q}(\O)} \le \cter{afdinsob} \norm{u}{1,p,\mesh}\qquad  \forall u \in H_\disc(\O),
\label{insob}
\ee
where $\dsp \norm{u}{1,p,\mesh}^p$ is defined in \refe{defnorunpd}.
\label{sob}
\end{lemma}
\begin{proof}
If $p=1$, one takes $q=1^\star$ and the result follows from Lemma~\ref{sobun} (in this case $\cter{afdinsob}$ does not depend on $\eta$).
If $1<p<d$, one takes $q=p^\star$ applies Lemma~\ref{sobp}.

\vspace{0.2cm}

If $p \ge d$, one chooses any $q \in ]p, \infty[$ and $p_1 < d$ such that $p_1^\star=q$ (this is possible since $p_1^\star$ tends to $\infty$ as $p_1$
tends to $d$). 
Lemma~\ref{sobp} gives, for some $\cter{afdinsobp}$ only depending on $p$, $d$ and $\eta$, that $\norm{u}{L^{q}(\O)} \le \cter{afdinsobp} \norm{u}{1,p_1,\mesh}$.
But, using H\"older's inequality, there exists $\ctel{afdc2}$, only depending on $d$, $p$, $\O$, such that $\norm{u}{1,p_1,\mesh} \le \cter{afdc2} \norm{u}{1,p,\mesh}$. 
Inequality \refe{sob} follows with $\cter{afdinsob}=\cter{afdinsobp} \cter{afdc2}$.
\end{proof}

\subsection{Compactness results for bounded families in the discrete $W^{1,p}$ norm}

\subsubsection{Compactness in $L^p$}

We prove in this section that bounded families in the discrete $W^{1,p}$ norms are relatively compact in $L^p$. 
We begin here also with the case $p=1$, giving in this case a crucial inequality which holds for general polyhedral partitions of $\O$.

\begin{lemma}Let $d \ge 1$ and let $\O$ be an open bounded set in $\R^d$, whose boundary is a finite union of parts of hyperplanes. 
Let $\mesh$ be a polyhedral partition of $\O$ in the sense of Definition \ref{defpart}. 
Then, with the notations of Definition \ref{defpart},
\begin{equation}
\norm{u(\cdot +\y) - u}{L^{1}(\R^d)} \le \vert \y \vert \sqrt{d} \norm{u}{1,1,\mesh} \qquad \forall u \in H_\mesh(\O), \, \, \forall  \y \in \R^d,
\label{inprekolm}
\end{equation}
where $u$ is defined on the whole space $\R^d$, taking $u=0$ outside $\O$, and $\vert h \vert$ is the Euclidean norm of $h \in \R^d$.
\label{prekolm}
\end{lemma}
\begin{proof}
One may prove this result in a similar way to that of \cite[Lemma 9.3 p. 770]{book} where an $L^2$ estimate on the translations is proven. 
Indeed, the proof of Lemma 9.3 \cite{book} holds  in the case  $p=1$ considered here for a general partition, while for  $p>1$, it requires the orthogonality condition satisfied by the admissible meshes of \cite[Definition 9.1 p 762]{book}.
We give here a simpler proof dedicated to the case $p=1$, using the $BV-$space, as in Lemma \ref{sobun}.

\vspace{0.2cm}

Let $u \in C^\infty_c(\R^d)$. For $\x, \y \in \R^d$, one has:
$$
\vert u(\x+\y)-u(\x) \vert = \vert \int_0^1 \grad u (\x+t\y) \cdot \y \d t \vert \le 
\vert \y \vert \int_0^1 \vert  \grad u (\x+t\y) \vert \d t.
$$
Integrating with respect to $\x$ and using Fubini's Theorem gives the  well-known  result
\be
\norm{u(\cdot+\y)-u}{L^{1}(\R^d)} \le \vert \y \vert  \int_{\R^d} \vert  \grad u \vert \d\x
\le \vert \y \vert \sum_{i=1}^d \norm{D_i u}{L^1(\R^d)},
\label{twuu}
\ee
where $\grad u=(D_1 u, \ldots, D_d u)$. By density of $C^\infty_c(\R^d)$
in $W^{1,1}(\R^d)$, Inequality~\refe{twuu} is also true
for $u \in W^{1,1}(\R^d)$.

\vspace{0.2cm}

We proceed now as in Lemma \ref{sobun}, using the same notations. 
Let $u \in BV$ and $u_n= u \star \rho_n$. 
Since $u_n \in W^{1,1}(\R^n)$, Inequality \refe{twuu} gives, for all $\y \in \R^d$, 
$\norm{u_n(\cdot+\y)-u_n}{L^{1}(\R^d)} \le \vert \y \vert  \sum_{i=1}^d \norm{D_i u_n}{L^1(\R^d)}$.
But, for $i=1,\ldots,d$, as in Lemma~\ref{sobun}, $\norm{D_i u_n}{L^1(\R^d)} \le \norm{D_iu}{M}$.
Then, since $u_n \tends u$ in $L^1(\R^d)$, as $n \tends \infty$, we obtain:
\be
\norm{u(\cdot+\y)-u}{L^{1}(\R^d)} \le \vert \y \vert  \sum_{i=1}^d \norm{D_i u}{M} =
\vert \y \vert \norm{u}{BV}\qquad \forall u \in BV, \, \, \forall \y \in \R^d.
\label{tbv}
\ee
Let $u \in H_\mesh(\O)$. One sets $u=0$ outside $\O$ so that $u \in L^1(\R^d)$; thanks to lemma \ref{sobun}, $\norm{u}{BV} \le \sqrt{d}\norm{u}{1,1,\mesh}$ and thus:
$$
\norm{u(\cdot+\y)-u}{L^{1}(\R^d)} \le \vert \y  \vert \sqrt{d} \norm{u}{1,1,\mesh}\qquad\forall \y \in \R^d.
$$
\end{proof}

An easy consequence of Lemmas \ref{sobun} and \ref{prekolm} is a compactness result in $L^1$ given in the following lemma.

\begin{lemma}Let $d \ge 1$ and let $\O$ be an open bounded set in $\R^d$, such that its boundary $\partial \Omega$ is a finite union of parts of hyperplanes. 
Let $\cal F$ be a family of  polyhedral partitions of $\O$ in the sense of Definition \ref{defpart}. 
For $\mesh \in {\cal F}$, let $u_\mesh \in H_\mesh(\O)$ and assume that there exists $C \in \R$ such that for all $\mesh \in {\cal F}$, $\norm{u_\mesh}{1,1,\mesh} \le C$.
Then, the family $(u_\mesh)_{\mesh \in {\cal F}}$ is relatively compact in $L^1(\O)$ and also in $L^1(\R^d)$ taking $u_\mesh=0$ outside $\O$.
\label{kolmun}
\end{lemma}
\begin{proof}
By Lemma \ref{sobun}, the family  $(u_\mesh)_{\mesh \in F}$ is bounded in $L^{1^\star}(\O)$. 
Since $\O$ is bounded, the family $(u_\mesh)_{\mesh \in F}$ is bounded in $L^1(\O)$ and also in $L^1(\R^d)$, taking $u_\mesh=0$ outside $\O$. 
Thanks to the Kolmogorov compactness theorem, Lemma~\ref{prekolm} gives that the  family $(u_\mesh)_{\mesh \in {\cal F}}$ is relatively compact in $L^1(\O)$ and also in $L^1(\R^d)$ taking $u_\mesh=0$ outside $\O$.
\end{proof}

Note that in fact, the above result also holds for general (non polyhedral) partitions of $\O$, for instance in the case of curved boundaries. 
In the case $p >1$, we need an additional hypothesis on the meshes which we state in the following lemma.

\begin{lemma}Let $d \ge  1$, $1\le p< \infty$ and $\O$ be a polyhedral open bounded connected subset of $\R^d$.
 Let $F$ be a family of meshes of $\O$ in the sense of Definition~\ref{adisc}.
Let $\eta >0$ be such that, for all $\disc \in F$, one has $\eta \le d_{K,\sigma}/d_{L,\sigma} \le 1/\eta$ for all $\sigma \in \edges$, where $\mesh_\sigma=\{K,L\}$. 
For $\disc \in F$, let $u_\disc \in H_\disc(\O)$ and assume that there exists $C \in \R$ such, for all $\disc \in F$, $\norm{u_\disc}{1,p,\mesh} \le C$.
Then, the family $(u_\disc)_{\disc \in F}$ is relatively compact in $L^p(\O)$ and also in $L^p(\R^d)$ taking $u_\disc=0$ outside $\O$.
\label{kolmp}
\end{lemma}
\begin{proof}
Thanks to Lemma~\ref{sob} and to the fact that $\O$ is bounded, the family $(u_\disc)_{\disc \in F}$ is bounded in $L^1(\O)$ and also in $L^1(\R^d)$ taking $u_\disc=0$ outside $\O$.
Thanks once again to the fact that $\O$ is bounded, the family $(\norm{u_\disc}{1,1,\mesh})_{\disc \in F}$ is bounded in $\R$.
Then, as in the previous lemma, the Kolmogorov compactness theorem gives that the family $(u_\disc)_{\disc \in F}$ is relatively compact in $L^1(\O)$ and also in $L^1(\R^d)$ taking $u_\disc=0$ outside $\O$.

\vspace{0.2cm}

In order to conclude we use, once again, Lemma~\ref{sob}. 
It gives that the family $(u_\disc)_{\disc \in F}$ is bounded in $L^q(\O)$ for some $q>p$. 
With the relative compactness in $L^1(\O)$, this leads to the fact that the family $(u_\disc)_{\disc \in F}$ is relatively compact in $L^p(\O)$ (and then also in $L^p(\R^d)$ taking $u_\disc=0$ outside $\O$).
\end{proof}

\subsubsection{Regularity of the limit}

With the hypotheses of Lemma \ref{kolmp}, assume that $u_\disc \tends u$ in $L^p$ as $\size(\disc) \tends 0$ (Lemma~\ref{kolmp} gives that this is possible, at least for subsequences of sequences of meshes
with vanishing size). 
We prove below that $u \in W^{1,p}_0(\O)$.

\begin{lemma}Let $d \ge  1$, $1\le p< \infty$ and let $\O$ be a polyhedral open bounded connected subset of $\R^d$.
Let $(\disc_n)_\nnn$ be a family of discretisations of $\O$ in the sense of Definition~\ref{adisc}.
Let $\eta >0$ be such that, for any discretisation $\disc_n = (\mesh_n,\edges_n,\points_n)$,  one has $\eta \le d_{K,\sigma}/d_{L,\sigma} \le 1/\eta$ for all
$\sigma \in \edges$, where $\mesh_\sigma=\{K,L\}$. 
For $n \in \N$, let $u^{(n)}  \in H_{\disc_n}(\O)$ and assume that there exists $C \in \R$ such, for all $n \in \N$, $\norm{u^{(n)}}{1,p,\mesh_n} \le C$.
Assume also that $\size(\disc_n) \tends 0$ as $n \tends \infty$. 
Then:
\begin{enumerate}
\item There exists a sub-sequence of $(u^{(n)})_\nnn$, still denoted by $(u^{(n)})_\nnn$, and $u \in L^p(\O)$ such that $u^{(n)} \tends u$ in $L^p(\O)$ as $n \tends \infty$.
\item $u \in W^{1,p}_0(\O)$ and 
\be
\norm{\grad u }{L^p( \O)^d} = \norm{\ \vert \grad u \vert \ }{L^p(\O)} \le \frac {(1+\eta) d^{\frac {p-1} p}} {\eta} C
\label{majnormlimit}\ee
(recall that $\vert  \grad u \vert$ is the Euclidean norm of $\grad u$).
\end{enumerate}
\label{reglim}
\end{lemma}
\begin{proof}
The fact that there exists a subsequence of $(u^{(n)})_\nnn$, still denoted by $(u^{(n)})_\nnn$, and $u \in L^p(\O)$ such that $u^{(n)} \tends u$ in $L^p(\O)$ as $n \tends \infty$ is a consequence of the relative compactness of $(u^{(n)})_\nnn$ in $L^p$ given in Lemma~\ref{kolmp}. 
There only remains to prove that $u \in  W^{1,p}_0(\O)$.

\vspace{0.2cm}

Letting  $u^{(n)}=0$ and $u=0$ outside $\O$, one also has $u^{(n)} \tends u$ in $L^p(\R^d)$.
Let us now construct an approximate gradient, denoted by $\tilde{\grad}_{\disc_n} u^{(n)}$, bounded in $L^p(\O)$, equal to $0$ outside $\O$
and converging, at least in the distributional sense, to $\grad u$.

\vspace{0.2cm}

{\bf Step 1} Construction of $\tilde{\grad}_\disc u$, for $u \in H_\disc(\O)$, and its properties.

Let $\nnn$ and $\disc=\disc_n$. For this step, one sets $u=u^{(n)}$ (not to be
confused with the limit of the sequence $(u^{(n)})_\nnn$).
For $\sigma \in \edges$, one sets $u_\sigma=0$ if $\sigma$ is on the boundary
of $\O$. 
Otherwise, one has $\mesh_\sigma = \{K,L\}$ and we choose  a value $u_\sigma$ between $u_K$ and $u_L$  (it is possible to choose, for instance, $u_\sigma = \frac 1 2 (u_K + u_L)$ but any other choice between $u_K$ and $u_L$ is possible). 
Then, one defines $\tilde{\grad}_\disc u$
on $K \in \disc$ in the following way:
$$
\tilde{\grad}_\disc u=\frac 1 {\mcv} \sum_{\sigma \in \edgescv}
\medge \ncvedge (u_\sigma - u_K).
$$
The function $\tilde{\grad}_\disc u$ is constant on each $K \in \mesh$ and, on $K$, using H\"older's inequality
$$
\vert \tilde{\grad}_\disc u \vert^p \le \frac 1 {(\mcv)^p} \left(\sum_{\sigma \in \edgescv} \medge \ncvedge |u_\sigma - u_K|\right)^p \le \frac 1 {(\mcv)^p}
\left(\sum_{\sigma \in \edgescv} \medge d_{K,\sigma}\right)^{p-1}
\sum_{\sigma \in \edgescv} \medge d_{K,\sigma}\left(\frac {D_\sigma u} {d_{K,\sigma}}\right)^p.
$$
Since $\sum_{\sigma \in \edgescv} \medge d_{K,\sigma} =d \mcv$, one deduces
$$
\vert \tilde{\grad}_\disc u \vert^p \le  \frac {d^{p-1}} {\mcv}\sum_{\sigma \in \edgescv} \medge d_{K,\sigma} \left(\frac {D_\sigma u} {d_{K,\sigma}}\right)^p.
$$
This gives an $L^p$- estimate on $\tilde{\grad}_\disc u$ in $(L^p(\O))^d$ (or in $(L^p(\R^d))^d$, setting $\tilde{\grad}_\disc u=0$ outside $\O$), in terms
of $\norm{u}{1,p,\mesh}$, namely
\be
\norm{\vert \tilde{\grad}_\disc u \vert}{L^p} \le \frac {(1+\eta) d^{\frac {p-1} p}} {\eta} \norm{u}{1,p,\mesh}.
\label{estgrad}
\ee
In order to prove, in the next step, the convergence of this approximate gradient, we now compute the integral of this gradient against a test function.
Let $\varphi \in C^\infty_c(\R^d;\R^d)$, $\varphi_K$ the mean value of $\varphi$ on $K \in \disc$, and $\varphi_\sigma$ the mean value of $\varphi$ on $\sigma$. Then,
\be
\int _{\R^d} \tilde{\grad}_\disc u \cdot \varphi \d\x= \sum_{K \in \disc}
\sum_{\sigma \in \edgescv}
\medge \ncvedge (u_\sigma - u_K) \varphi_K=
\sum_{K \in \disc} \sum_{\sigma \in \edgescv}
\medge \ncvedge (-u_K) \varphi_\sigma + R(u,\varphi),
\label{intcphi}
\ee
with
$$
R(u,\varphi)=\sum_{K \in \disc}
\sum_{\sigma \in \edgescv}
\medge \ncvedge (u_\sigma - u_K) (\varphi_K-\varphi_\sigma).
$$
Then, there exists $C_\varphi$ only depending on $\varphi$, $d$, $p$, $\O$ and $\eta$ such that $\vert R(u,\varphi) \vert \le C_\varphi \size(\disc) \norm{u}{1,p,\mesh}$. 
Equation \refe{intcphi} can also be written as

\be
\int _{\R^d} \tilde{\grad}_\disc u \cdot \varphi \d\x=
\sum_{K \in \disc} \int_K (-u_K)\  \div(\varphi)\  \d\x + R(u,\varphi)= -\int_{\R^d} u\  \div(\varphi)\  \d\x
+ R(u,\varphi).
\label{intcphid}
\ee

\vspace{0.2cm}

{\bf Step 2} Convergence of $\tilde{\grad}_{\disc_n} u^{(n)}$ to $\grad u$ and proof of
$u \in W^{1,p}_0(\O)$ .

We consider now the sequence $(u^{(n)})_\nnn$. Inequality \refe{estgrad} gives
$$
\norm{\vert \tilde{\grad}_\disc u^{(n)} \vert}{L^p} \le \frac {(1+\eta) d^{\frac {p-1} p}} {\eta}
 \norm{u^{(n)}}{1,p,\mesh}.
 $$
Then, the sequence $(\tilde{\grad}_\disc u^{(n)} )_\nnn$ is bounded in $L^p(\R^d)^d$ and we can assume, up to a subsequence, that $\tilde{\grad}_\disc u^{(n)}$ converges to some $w$ weakly in $L^p(\R^d)^d$, as $n \tends \infty$ and
$\norm{\ \vert w \vert \ }{L^p} \le \frac {(1+\eta) d^{\frac {p-1} p}} {\eta} C$.

\vspace{0.2cm}

Let $\varphi \in C^\infty_c(\R^d;\R^d)$, Equation \refe{intcphid} gives
\be
\int _{\R^d} \tilde{\grad}_\disc u^{(n)} \cdot \varphi \d\x= -\int_{\R^d} u^{(n)} \ \div(\varphi) \ \d\x
+ R(u^{(n)},\varphi).
\label{graddiv}
\ee
Thanks to $\vert R(u^{(n)},\varphi) \vert \le C_\varphi \size(\disc_n) \norm{u^{(n)}}{1,p,\mesh_n}$,
one has $R(u^{(n)},\varphi) \tends 0$, as $n \tends \infty$.
Since $u^{(n)} \tends u$ in $L^p(\R^d)$ as $n \tends \infty$, passing to the limit
in \refe{graddiv} gives:
$$
\int _{\R^d} w \cdot \varphi \d\x= -\int_{\R^d} u \ \div(\varphi) \ \d\x.
$$
Since $\varphi$ is arbitrary in $C^\infty_c(\R^d,\R^d)$, one deduces
that $\grad u = w$. Then $u \in W^{1,p}(\R^d)$ and
$\norm{\vert \grad u \vert}{L^p} \le \frac {(1+\eta) d^{\frac {p-1} p}} {\eta} C$.
Finally, since $u=0$ outside $\O$, one has $u \in W^{1,p}_0(\O)$.
\end{proof}

\section{Conclusion and perspectives}\label{conc}
A symmetric discretisation scheme was introduced for anisotropic heterogeneous problems on distorted nonconforming meshes. Although this scheme stems from the finite volume analysis, which was developed these past years, its formulation is actually derived from a discrete weak formulation; in this respect it may be seen as a nonconforming finite element method. 
Tools of functional analysis were obtained, which allow a mathematical analysis of the scheme; the convergence of the discrete solution to the exact solution of the continuous problem is shown with no regularity assumption on the solution (other than the natural assumption that it is in  $H^1_0(\Omega)$). 
Even though this convergence result yields no rate of convergence, it is probably more interesting than error estimates which require some assumptions on the diffusion tensor. 
Nevertheless, we show an order 1 estimate in the case of the Laplace operator, which is readily adaptable to regular (say piece-wise $C^1$) isotropic diffusion operators.  
The numerical results presented here show the good performance of the scheme (in particular order 2 is obtained for the convergence in the $L^2$ norm of the solution), and so do three dimensional experiments which were performed in \cite{che-09-col} for the incompressible Navier--Stokes equations on general grids.
Note that the convergence analysis which is performed here readily extends to the non-linear setting of Leray-Lions operators. This will be the subject of a future paper.

\end{document}